%% file: records.tex
\documentclass[11pt]{amsart}

\usepackage{epsfig,amsmath}

\theoremstyle{definition}

\def\dm#1#2#3 {{\em Discrete Math.\ }{\bf#1} (#2) #3}
\def\jspi#1#2#3 {{\em J. Statist.\ Plann.\ Inference }{\bf#1} (#2) #3}

\def\md{\mathrm{md}}
\def\ub{\mathrm{ub}}

\numberwithin{equation}{section}

\begin{document}

\parskip 6pt
\parindent 0pt
\baselineskip 13pt

\title[Lower bounds for maximal determinants]{New Lower Bounds
for the Maximal Determinant Problem}

\author[Orrick]{William P. Orrick}
\address{Department of Mathematics, Indiana University,
Bloomington, IN 47405}
\email{worrick@indiana.edu}

\author[Solomon]{Bruce Solomon}
\address{Department of Mathematics, Indiana University,
Bloomington, IN 47405}
\email{solomon@indiana.edu}

\author[Dowdeswell]{Roland  Dowdeswell}
\address{11 West 17th St., Floor 2, New York, NY 10011}
\email{elric@imrryr.org}

\author[Smith]{Warren D. Smith} \address{Department of
Mathematics, Temple University, Philadelphia, PA 19122}
\email{wds@math.temple.edu}

\date{First draft April 7, 2003. Last Typeset \today.}

\begin{abstract}
We report new world records for the maximal determinant of an
$n\times n$ matrix with entries $\,\pm 1\,$. Using various
techniques, we beat existing records for $n=22$, $23$, $27$,
$29$, $31$, $33$, $34$, $35$, $39$, $45$, $47$, $53$, $63$,
$69$, $73$, $77$, $79$, $93$, and $95$, and we present the
record-breaking matrices here. We conjecture that our $n=22$
value attains the globally maximizing determinant in its
dimension. We also tabulate new records for $n=67$, $75$, $83$,
$87$, $91$ and $99$, dimensions for which no previous claims
have been made. The relevant matrices in all these dimensions,
along with other pertinent information, are posted at
\texttt{www.indiana.edu/$\sim$maxdet}.
\end{abstract}

\maketitle

\section{Introduction}\label{sec:intro}

Let $\,\md(n)\,$ denote the global maximum of the determinant
function on the set of all $\,n\times n\,$ matrices populated
entirely by $\,\pm1$s. It is easy to show that $\,\md(n)\le
n^{n/2}\,$, and the well-known Hadamard conjecture predicts that
when $\,n>2\,$, equality occurs if $\,n\equiv 0\pmod 4\,$.
This conjecture has been confirmed for all $\,n\le 424\,$ and
many larger $\,n\,$, including a number of infinite
sequences~\cite{sy}.

When $\,n>2\,$ is \emph{not} divisible by 4, however, the
inequality is necessarily strict, and relatively few exact
values for $\,\md(n)\,$ are known.

Using a combination of techniques, we find a number of
matrices whose determinants exceed previous records, thereby
establishing new lower bounds for $\,\md(n)\,$. We sketch our
methods Sec.\ \ref{sec:methods} below, and will expand on them
further in a forthcoming survey paper~\cite{OOS}. Here, our
purpose is simply to record and establish priority for our new
lower bounds.

Other authors have derived \emph{upper} bounds for $\,\md(n)\,$.
When $\,n\not\equiv 0\pmod 4\,$ these are more stringent than
the simple Hadamard estimate $\,\md(n)\le n^{n/2}\,$. We compare
our lower bounds to these upper bounds below, using formulae
which come from Barba~\cite{Ba} when $\,n\equiv1\pmod 4\,$, from
Ehlich~\cite{Eh1} and Wojtas~\cite{Wo} (independently) when
$\,n\equiv 2\pmod 4\,$, and from Ehlich~\cite{Eh2} when
$\,n\equiv 3\pmod 4\,$.

These bounds cannot be sharp in the cases we treat here, except
possibly when $\,n=91\,$. Indeed, we conjecture that
$\,\md(22)=2^{23}\times5^{11}\,$, i.e. that the value we present
for $\,n=22\,$ below is in fact best possible, though it achieves
only 95\% of the theoretical bound.

The upper bound formulae we use, as well as the matrices
associated with all determinants listed in Sec.~\ref{sec:records}
below, are available at
\texttt{http://www.indiana.edu/$\sim$maxdet}, a website we have
created as an archive for data pertaining to the determinant
maximization problem. We intend to update it with increasingly
better values, matrices, and references as they become
available.

After a brief summary of our methods briefly in
Sec.\ \ref{sec:methods} below, we present some of our findings in
Sec.\ \ref{sec:records}. Our records in dimensions $n=22$,
$23$, $27$, $29$, $31$, $33$, $34$, $35$, $39$, $45$, $47$,
$53$, $63$, $69$, $73$, $77$, $79$, $93$, and $95$ either
eclipse a previous record, or manifest visible structure, or
both. Except for the $n=93$ case (which barely misses) each of
these
examples also exceeds 75\% of the theoretical upper bound in its
dimension. Many do much better. In these cases we therefore
present our record-breaking matrix of $\,\pm 1$s, along with its
determinant, a comparison to the theoretical bound, and a brief
note, including a reference to the previous record.

Our work has also produced matrices of large determinant in
dimensions not listed above. While they may not beat any
previously published records, we believe several of them do
provide best-known lower bounds for $\,\md(n)\,$. We list such
values for $\,n=67$, $75$, $83$, $87$, $91$ and $99$ in a table
at the end of Sec. \ref{sec:records}.

One further note: In dimensions $n=29$, $33$, $45$, $53$, $69$,
$73$, $77$ and $93$ (each congruent to 1 mod 4), our
determinants improve on values published by Farmakis \& Kounias,
who constructed their matrices by adding a suitable row and
column to a Hadamard matrix of largest known excess \cite{fk}.
The \emph{excess} of a Hadamard matrix is simply the sum of its
entries, and has been studied in connection with the
classification problem for Hadamard matrices. By deleting the
appropriate row and column from matrices we have found in
dimensions $73$ and $77$, one gets Hadamard matrices with larger
excess than any previously known in those dimensions. For
instance, we get a $72\times72$ Hadamard with excess 580,
improving the 576 given in \cite{fk}, by deleting the first row
and column from the $n=73$ example in Sec.\ \ref{sec:records}
below. Proceeding similarly with an $n=77$ example posted on our
website, one gets a $76\times76$ Hadamard with excess 628,
beating the value of 620 given in \cite{fk}. We present a
different $n=77$ example in Sec.\ \ref{sec:records} because its
determinant is larger.

\section{Methods.}\label{sec:methods}

Where not otherwise noted, the matrices in Sec.~\ref{sec:records}
were discovered numerically, using a hill-climbing computer
program which combines a discrete version of gradient-ascent
with simulated annealing. The first program of this type was
used by Smith, and described in his thesis~\cite{Sm}.

Stated roughly, the algorithm runs as follows. First one chooses
a starting matrix; a random array of $\,\pm 1$s, or better, a
well-chosen guess as to what the maximizer might look like. One
then passes repeatedly through this candidate, modifying it row
by row to raise its determinant. Since changes in a row don't
affect the cofactors of that row, one can exploit the cofactor
expansion formula for determinants to do this very efficiently.

This simple strategy encounters a fundamental limitation,
however: It almost always stalls in a \emph{local} maximum of
the determinant function, far below $\,\md(n)\,$. We over- come
this problem to some extent by randomly perturbing the ascent
when it gets stuck in this way, and attempting to continue from
there.

By carefully examining our best result found using this method in
the case $\,n=31\,$, we discovered an explicit construction which
produces new records for $\,n=47\,$, $63$, $79$ and $95$. More
generally, the construction applies when $\,n\equiv 15\pmod
{16}\,$, and uses conference matrices~\cite{ms}. We describe it
briefly as follows.

Suppose $\,n\equiv 15\pmod {16}\,$, and write $\,n=4k+3\,$. Then
$\,k+1\equiv 0\pmod 4\,$, and (at least in the dimensions of
interest here) there exists an antisymmetric
$\,(k+1)\times(k+1)\,$ conference matrix $\,C=[c_{ij}]\,$.
Normalize $\,C\,$, by multiplying rows (and corresponding
columns, to preserve antisymmetry) by $\,-1\,$ so that
$\,c_{1j}=1\,$ for all $\,j>1\,$, and $\,c_{i1}=-1\,$ for all
$\,i>1\,$. Now delete the first row and column of $\,C\,$ to
produce a  $\,k\times k$ matrix $B$, with the property that
$BB^{T}= (k-1)\,I_{k-1} - J_{k-1}\,$, $\,J\,$ here denoting a
matrix comprised entirely of $\,+1$s.

Given $\,B\,$, we now obtain a $\,\pm 1\,$ matrix in four steps:
First, tensor $\,B\,$ with a $\,4\times 4\,$ Hadamard
matrix $\,H_4\,$. Second, replace each $0$ in the tensor product
by a $-1$. Third, pad this $\,4k\times4k\,$ matrix with 3
initial rows, obtained by normalizing $\,H_4\,$, deleting a row,
and then tensoring with a row of $\,k\,$ 1s. Finally, we
similarly pad with three initial columns. The examples in
Sec.~\ref{sec:records} should clarify this.

Note that if $\,k\,$ is a prime power, the Jacobsthal matrix of
the finite field $\,\mathrm{GF}(k)\,$ provides a matrix $\,B\,$
with exactly the properties we need.

In any case, this construction produces a matrix of very high
determinant when $\,n\,$ is not too big. As $\,n\to\infty\,$,
however, the resulting determinant tends to $0$ relative to
Ehlich's theoretical upper bound.
\vfill
\goodbreak
\section{Record-holding matrices.}\label{sec:records}

\input{matData.tex}


\textbf{Other cases of interest:} The table below gives the
largest determinants we have found among $\,\pm 1\,$ matrices in
several additional dimensions, with $\,\ub(n)\,$ denoting the
theoretical upper bound for $\,\md(n)\,$ mentioned in Sec.
\ref{sec:intro}. They do not quite achieve 75\% of the theoretical bound,
and to the best of our knowledge, no previous bounds have been
published. Further details, including the matrices themselves,
can be found at our website,
\texttt{www.indiana.edu/$\sim$maxdet}$\ $.
\smallskip

\begin{equation*}
\begin{array}[b]{cll}
\mathbf{Size} & \textbf{Lower bound for}\ \mathbf{md(n)}
& \textbf{Fraction of}\ \mathbf{ub(n)}\\[2pt] \\
n=67: &\ \quad  2^{66}\times 16^{31} \times 4765 &\qquad0.7677\\
n=75: &\ \quad  2^{74}\times 18^{35} \times 6064 &\qquad0.7303\\
n=83: &\ \quad  2^{82}\times 20^{38} \times 157788 &\qquad 0.7322\\
n=87: &\ \quad  2^{86}\times 21^{41} \times 8777 &\qquad 0.7220\\
n=91: &\ \quad  2^{90}\times 22^{43} \times 9826 &\qquad 0.7203\\
n=99: &\ \quad  2^{98}\times 24^{47} \times 12118 &\qquad 0.7160\\
\end{array}
\end{equation*}
\bigskip

\section*{Acknowledgements}
The websites of N.J.A. Sloane, J. Seberry, C. Koukouvinos, and
E. Spence provided us with Hadamard and other large-determinant
matrices that we used to construct starting points for our own
searches. We gratefully acknowledge them for this contribution.
We made extensive use of \emph{Mathematica} during this project,
and we also thank Indiana University for the use of their Sun
E10000 computing platform.  Finally, we thank Michael Neubauer
for helpful remarks and interest in our work.

\bigskip

\end{document}

%% file: matData.tex
\parskip 4pt
\parindent 0pt
\baselineskip 12pt

\vfill
$\,\mathbf{n=22}\,$:

\begin{minipage}[t]{2in}
\baselineskip 0pt
\begin{verbatim}
++-----++++++++++-----
++-----+++++-----+++++
--+----+----++-----+-+
---+----+----++--+--+-
----+----+--+---+-+-+-
-----+----+---++--+--+
------+----+---+++-+--
+++----+------+-+++---
++-+----+---+--+--++--
++--+----+---+-+-+---+
++---+----+--+--+--++-
++----+----++-+-----++
+---+++++-+-+----+----
+-+--++-++-+-+----+---
+--+++-+-+-+--+----+--
+-++--++-++----+----+-
+-+++---+-++----+----+
-+++-+---++++----+----
-+-++-++--++-+----+---
-++-+-+-+++---+----+--
-++-++-++--+---+----+-
-+-+-+++++------+----+
\end{verbatim}
\end{minipage}
\
\makebox[0.5in]{}
\
\begin{minipage}[t]{3in}
\begin{eqnarray*}
\text{Determinant:}&&2^{23}\times5^{11}\\
\text{Theoretical bound:}&&2^{21}\times 3\times5^{10}\times 7\\
\text{Fraction of bound:}&&0.95
\end{eqnarray*}
\vskip 0.21in

We conjecture that the determinant above \emph{is} $\,\md(22)\,$.
 Smith~\cite{Sm}, and more recently, Cohn~\cite{co1,co2}
published earlier records. Cohn's value was
 $\,2^{21}\times3^2\times23^2\times197^2\,$ (90.1\% of bound).
\vskip 0.05in
\hrule
\end{minipage}

\vfill
$\,\mathbf{n=23}\,$:

\begin{minipage}[t]{2in}
\baselineskip 0pt
\begin{verbatim}
+-+----++----++++++++++
-++--++----++--++++++++
++-++----++----++++++++
--++---++++++--++--++--
--+-+--++++++----++--++
-+-----+++------++--++-
-+-----++-+----+--++--+
+--++++++--+---+--+-++-
+--++++++---+---++-+--+
--++++----------+-++-+-
--+++-+--------+-+--+-+
-+-++--+---+++++-+-+-+-
-+-++---+--++++-+-+-+-+
+----------+++---++++--
+----------++-+++----++
++++--++--++--+-+++----
++++-+--++--+-++-++----
+++-++--+-++-+-+++-----
+++-+-++-+--++-++-+----
++++--+-++-+-+-----+-++
++++-+-+--+-++------+++
+++-++-+-+-+--+----++-+
+++-+-+-+-+-+-+----+++-
\end{verbatim}
\end{minipage}
\
\makebox[0.5in]{}
\
\begin{minipage}[t]{3in}
\begin{eqnarray*}
\text{Determinant:}&&2^{22}\times3\times 5^6\times67\times211\\
\text{Theoretical bound:}&&2^{22}\times3\times5^6\times675\sqrt{505}\\
\text{Fraction of bound:}&&0.931983
\end{eqnarray*}
\vskip 0.85in

Smith~\cite{Sm} published an earlier record of
$\,2^{22}\times19\times5741^2\,$ (88\% of bound).
\vskip 0.05in
\hrule
\end{minipage}
\vfill
\eject

\phantom{.}
\vfill
$\,\mathbf{n=27}\,$:

\begin{minipage}[t]{2in}
\baselineskip 0pt
\begin{verbatim}
-++--+-++--+-++++++++-+-++-
++--+--+-+--++++++++++---++
+--+-+---+++-++++++++-++--+
--+++---+-+-++++++++++-++--
-+-+--++--+---++++++----+-+
+-+---+-++----++++++---+-+-
----++-++++---++++--+++----
-++-+-++++++++-++--+--++--+
+-+--+++++++++++--+--+--+-+
-+-+-++++++++++--+-+-+-+-+-
+--++-++++++++--+++---+-++-
--++---++++++--+-+--+----++
++-----++++++-+-+---+--++--
++++---++++--+----++-++----
+++++++-++--+-+--+-++-+-+-+
+++++++++--+----+++-++-+--+
++++++++--+-+-++--+-+-++-+-
+++++++--+++---++--+++--++-
++++++--+-+--++-+--------++
++++++-+-+---+-+-+-----++--
++++--+----++-++++---++----
++----+-+-+--+-+-+--+++++++
--++--++-+---++-+---+++++++
+-+--+-+--+-+---++-+-++++++
+--++--++--+--++---+-++++++
-+-+-+--++--+--++-+--++++++
-++-+----+++--+--++--++++++
\end{verbatim}
\end{minipage}
\
\makebox[0.5in]{}
\
\begin{minipage}[t]{3in}
\begin{eqnarray*}
\text{Determinant:}&&2^{26}\times6^{11}\times 518\\
\text{Theoretical bound:}&&2^{26}\times 5\times6^{10}\times44\sqrt{237}\\
\text{Fraction of bound:}&&0.917665
\end{eqnarray*}
\vskip 1.225in

Smith~\cite{Sm} published an earlier record of
$\,2^{27}\times11\times90481^2\,$ (88\% of bound).
\vskip 0.05in
\hrule
\end{minipage}

\vfill
$\,\mathbf{n=29}\,$:

\begin{minipage}[t]{2in}
\baselineskip 0pt
\begin{verbatim}
+---+++++++------++++++++++++
-+--+++++++++++++------++++++
--+-+++++++++++++++++++------
-----------++++++++++++++++++
+++-+--++++--++++--++++--++++
+++--+-++++++--++++--++++--++
+++---+++++++++--++++--++++--
+++-++++---+-+-+-++--+---+++-
+++-+++-+--+--+-++-+--++-+--+
+++-+++--+--+-++---+++-++--+-
+++-+++---+-++--+-+-+-+-+-+-+
-+++-++++--+++-----++++-++-++
-+++-++--+++--++--+-++++-++-+
-++++-+-++-+---++++++---+-+++
-++++-++--+-+--+++-+-++++++--
-+++++--+-+-+++--+++-+-+--+++
-+++++-+-+---++-+++-+-++++-+-
+-++-++++----+++++--++-++-+-+
+-++-++-++--++-++-++--++-+++-
+-+++-++--++-+++--++--+++--++
+-+++-+--+++++--++--++-+-+-++
+-++++--+-+++-++-+--+-+-++++-
+-++++-+-+-++-+-+-++-+--+++-+
++-+-+++--+-+-+++++++----+-++
++-+-++--+++-++-++-+-++-+-++-
++-++-+++--++-+-+-+-++++--++-
++-++-+-++--++++-++--++-++--+
++-+++--+-++-+-++-++++-+++---
++-+++-+-+-+++-+-+-++-++--+-+
\end{verbatim}
\end{minipage}
\
\makebox[0.5in]{}
\
\begin{minipage}[t]{3in}
\begin{eqnarray*}
\text{Determinant:}&&2^{28}\times7^{12}\times320\\
\text{Theoretical bound:}&&2^{28}\times7^{12}\times49\sqrt{57}\\
\text{Fraction of bound:}&&0.865001
\end{eqnarray*}
\vskip 1.05in

Smith~\cite{Sm}, Farmakis \& Kounias~\cite{fk}, and most
recently, Koukouvinos~\cite{ko}, have published earlier records.
Koukouvinos' value is $\,2^{28}\times7^{13}\times43\,$
(81.4\% of bound).
\vskip .05in
\hrule
\end{minipage}

\vfill
\eject

\phantom{.}
\vfill
$\,\mathbf{n=31}\,$:

\begin{minipage}[t]{2.5in}
\baselineskip 0pt
\begin{verbatim}
+++--++--++--++--++--++--++--++
+++-+-+-+-+-+-+-+-+-+-+-+-+-+-+
++++--++--++--++--++--++--++--+
--+----+++-+++----++++----+---+
-+-----++-+++-+--+-++-+--+---+-
+------+-+++-++-+--+-++-+---+--
+++-----+++-++++----++++---+---
--+---+----+++-+++----++++----+
-+---+-----++-+++-+--+-++-+--+-
+---+------+-+++-++-+--+-++-+--
++++--------+++-++++----++++---
--+---+---+----+++-+++----++++-
-+---+---+-----++-+++-+--+-++-+
+---+---+------+-+++-++-+--+-++
++++---+--------+++-++++----+++
--++++----+---+----+++-+++----+
-+-++-+--+---+-----++-+++-+--+-
+--+-++-+---+------+-+++-++-+--
+++-++++---+--------+++-++++---
--+---++++----+---+----+++-+++-
-+---+-++-+--+---+-----++-+++-+
+---+--+-++-+---+------+-+++-++
++++----++++---+--------+++-+++
--++++----++++----+---+----+++-
-+-++-+--+-++-+--+---+-----++-+
+--+-++-+--+-++-+---+------+-++
+++-++++----++++---+--------+++
--++++-+++----++++----+---+----
-+-++-+++-+--+-++-+--+---+-----
+--+-+++-++-+--+-++-+---+------
+++-+++-++++----++++---+-------
\end{verbatim}
\end{minipage}
\
\makebox[0.3in]{}
\
\begin{minipage}[t]{3in}
\begin{eqnarray*}
\text{Determinant:}&&2^{30}\times7^{12}\times5324\\
\text{Theoretical bound:}&&2^{30}\times7^{12}\times144\sqrt{1589}\\
\text{Fraction of bound:}&&0.927499
\end{eqnarray*}

\vskip 1.075in
Though discovered numerically, our confer\- ence-matrix
method also produces this example (see Sec.~\ref{sec:methods}).
Smith~\cite{Sm} published an earlier record of $\,2^{30}\times
5^4\times7^2\times11^2\times29^2\times149^2\,$ (87\% of bound).
\vskip 0.05in
\hrule
\end{minipage}

\vfill
$\,\mathbf{n=33}\,$:

\begin{minipage}[t]{2.5in}
\baselineskip 0pt
\begin{verbatim}
-+-+++++++----+++++++++++++++++++
++--------+++++++++++++++++++++++
--+--+++++++++++--+++-++-+-++-+++
+-++++++---+-++--+++-+--+--++++++
+-+++++-+-+-+-++-+--+-+-++++-+-++
+-+++-++-++-+---++++++-+-++---+++
+-+++--+++-+-+-+++----++-+++++++-
+--++++-++-++++++--+++--+++-+-+--
+--+++-++++-++--+-+-+++++--+++--+
+-+--+++++++--++++++-++++-+--+---
-++--+++-+++++--++--++--++++++-+-
-+++++---+++++-++-++--+-+-+--++++
-++++-+-+-+++++++++--+-+--+++---+
-++++--++-+++++--+-+++++++---++--
++-+-++++-++-+-++++-+----+---++-+
+++-++-++-+-+-+-+-++-----++++++--
++++-++--+--+++-++-+--++-+--++--+
+++-+-+-++-+-+----+++--+++++-+--+
++++-++-+--++---+-+--+++++-+--++-
++--++-+++-++++--++--++--++----++
++--+++-+++-++-+-+-+-+-+---+-+++-
+++-+-+++---++-++++++-+-+---+--+-
++++--+-+++--++---+-+++---+-++++-
++++-+-+-+--++++-++-+--++-++--+--
++-+++++--++-++-+--++-++--++---+-
++++---++++--++++--+-+--++-+---++
++++-+-++--++--+---+++-+--+-++-++
++-+--+++++++----+-+--+-+-+++-+-+
++++++---+++---+-++++++--+-++----
++-++-++-++++-++--+----+++--++-+-
+++-++--++++--+-++--+--++---+-+++
+++-+-++-+-++-+++---+++----+-++-+
+++-++++--+--+-+-----++++++-+-+-+
\end{verbatim}
\end{minipage}
\
\makebox[0.3in]{}
\
\begin{minipage}[t]{3in}
\begin{eqnarray*}
\text{Determinant:}&&2^{32}\times 8^{14}\times441\\
\text{Theoretical bound:}&&2^{32}\times 8^{14}\times64\sqrt{65}\\
\text{Fraction of bound:}&&0.854677
\end{eqnarray*}

\vskip 1.8in

Farmakis \& Kounias~\cite{fk} published an earlier record of
$\,2^{32}\times8^{15}\times51\,$ (79\% of bound).
\vskip 0.05in
\hrule
\end{minipage}

\vfill
\eject

\phantom{.}
\vfill
$\,\mathbf{n=34}\,$:

\begin{minipage}[t]{2.5in}
\baselineskip 0pt
\begin{verbatim}
--++++++++++++++++++++++++++++++++
++----+++--+--+--++--+++++++++++++
++-------++-++-++--++-++++++++++++
++++++++--+-+---++-+-++-+--+-+++++
++++++-++-+--++++----+++-++++---++
++++++--++-++--+-++-+-++--++-++-++
+++++++--+-+-++---+++-+-++-++--+++
+++++++-+---++++-+-+---++++-++++--
++-++-+++--+++++++-++-+---------+-
+++--++-+++-++++-+++-+-----+-----+
+++-++++-+-+-+--++-+++-+--++++----
+++-++-++-++-+-++++-+---++-+--++--
+++++-++-++-+--+-++-++--+++-+---+-
++++-+-+++-++-+++---++--+---++-+-+
++-+++-++++--++---+++++-+-+--++---
++-+++--+++++---++++-+++-+--+--+--
++++-+++--+++-+-+-+++--+-++---+--+
+++++-+--+++-++++-+--+-+-----++++-
-++-+-++++++---+---+--+---+-+-++-+
-+-++--+-++++++--+-------+-++++--+
-+++--++++--++--+-+---+++--++-+---
-+--++++++--++--+-+------++--+-+++
-+-+-+++++++---+---+---+++-+-+--+-
+-+-+-+++++-+-+-----+-++-+-+-+-+--
+---++++-+-+++-+-----+++++----+--+
+--+-++++++--+---+--+--+----+-++++
-++--++-+-++++------+++--+--+++-+-
+-++---+++-+++-----+-+---+++--+++-
+--+-+++--++++-+--+---+---++++-+--
+-+-+--++-++++----++---++---++--++
-+-++-+-+-++++------++-++-++---+-+
+-++--+-++++-+--++----+-+++--+---+
+---+++-+++++-+-+-------+-+++-+-+-
-++--+-+-++++++--+----+++-+----++-
\end{verbatim}
\end{minipage}
\
\makebox[0.3in]{}
\
\begin{minipage}[t]{3in}
\begin{eqnarray*}
\text{Determinant:}&&2^{33}\times 8^{16}\times28\\
\text{Theoretical bound:}&&2^{33}\times 8^{16}\times33\\
\text{Fraction of bound:}&&0.848485
\end{eqnarray*}

\vskip 1.4in

We find no prior records in the literature, though one produces
a $\,\pm 1\,$ matrix with determinant $\,2^{33}\times8^{16}\times25\,$
(75.76\% of bound) by tensoring an $\,n=17\,$ maximizer with a
$\,2\times2\,$ maximizer.
\vskip 0.05in
\hrule
\end{minipage}

\phantom{.}
\vfill
$\,\mathbf{n=35}\,$:

\begin{minipage}[t]{2.5in}
\baselineskip 0pt
\begin{verbatim}
---++----++++------++++----++++++++
-+++-+++++-+-+-+-+---++++--+---+++-
-++-++++++--+-+-+-++-+-++----++--++
+-+++++--+++-+--+-++-+---++----++++
++-++++--++-+-++-+---++--+++-++--+-
-++++-++--+++++--++----+----+-++-+-
-+++++--+-+++++++-------+--+-+--+-+
-++---++-+-++--++++-++---+-+--+-+-+
-++--+--++-++--+++++--+---+-++-+-+-
+++++--+++-++++----++++--++--------
+--++++---+++--++++++++++----------
+-++-+++++++-++-++-++------++++----
++--+++++++-++++--+-+-+-----+--++-+
--++-++--+-+++-+--+-+-+++++-+++---+
-+--+++--+-++-+-++-++--++++++--++--
-+-+-+-++-+-++--++-+-+++-++-+-+-+-+
--+-++-++-++--++--+++-++-+++--+-++-
-++++-+++-+-----+-+-+++-++++++-+---
------+++-++++++-+-+++--+++--+-+-++
++-+----+++-++-++-+++--+++-+--++-+-
+-+-+---++++--++++---+++++--+-++--+
++-+---+-+++--+-+++-+-+++-+--+--+++
+-+-+--+-++-++-+-+++-+-++-++++--+--
-++---+---+--++++--+++++-+--++--++-
-++--+----+--++--++++++-+-++--++--+
---++--+-+---+++++++--++-+-+-+-+--+
---++---++---++++++-++--+-+-+-+-++-
++-+-+-+---+--++--++-+--++--+++++--
+-+-++-+----++--++--+-+-++---+++++-
+-++-+-+----+-+++---++-+--+++--+-++
+-++--+-+---+-+--+++--+-++-++---+++
++--+-++---+-+-++--+--+-+-+++-+--++
++--++--+--+-+---++-++-+-+-+++---++
+++++-+-+------+-+-++--+--+--++++-+
+-----+-+--++++-+-+--+++--++-++++--
\end{verbatim}
\end{minipage}
\
\makebox[0.3in]{}
\
\begin{minipage}[t]{3in}
\begin{eqnarray*}
\text{Determinant:}&&2^{34}\times9^{12}\times23\times5167\\
\text{Theoretical bound:}&&2^{34}\times 8^{15}\times98\sqrt{134}\\
\text{Fraction of bound:}&&0.840907
\end{eqnarray*}

\vskip 2.15in

We find no prior records in the literature.
\vskip 0.05in
\hrule
\end{minipage}

\vfill
\eject

\phantom{.}
\vfill
$$\mathbf{n=39:}$$
\medskip

\begin{equation*}
\begin{minipage}[t]{3in}
\baselineskip 0pt
\begin{verbatim}
+-++--+---++++++++-+----+++-+----------
-+-++++-+--+---++-------+----+++--+--+-
-+++--++++--++-+--+---+-+-+---+--+--+++
-++++++++-+--+--++-+---++-++--++++-++--
--++++-+++++---+---+--+--+--+--+++--++-
--+++-+++------+-++--------++---+-++--+
-++-+++++---+++----++----+--+-+---++++-
-+++-++-++-+--+--+-+-+-+--+-----+-+-+++
-++++++--++---+--++-+---++-----+-+-+-++
++-+++++++++++--+-+++------+----+----++
++++++-++++--+++-+++-+-------++--++----
+++-++++++-++-+++++---++-------+---++--
--+--+-+--+-+--++++++--+--+--+-------+-
---++-+-+++-+++++++++-++++---+--+-+-+--
-+---+++---+-++++++++++-+---+---++-----
--++++--+--++++++++++++--+++--++------+
-+-+---++++++-++++++-+-+++-++-+----+-+-
----++--+-+++++--++---+++-+++----++--+-
-++-----+-++-+++--+++--++--+++-+----+-+
-+--+-+---++-+-+-+++-++--+++-+-----+++-
++++-++-+---+--+--++++++++++++-+++++-+-
++-+++-+------++-+--+-+++++++++-+---+++
+-++-+--++---+--++--+++-+--+++-----+++-
++-++---+-------+-+-++-+-++-+----+--+--
++-----+++-------+-++-+-++++---+--+----
+---+++-+++----+---++++++-+-+-+----+--+
+-++-+++--++------+-++++++-+--+---+-+--
+-+-+-+++-----+-+--+-+++++-+-+---+---++
+-+-+-++++-+-+---++-++-++++-+++++----+-
+--+---------++---++--++------+++--+-+-
+++-----+-++-+-+++--+-++-+----+-++++-++
+-+-+----+-+--+++-+++---+-++--+-++++++-
+--+--+++-++--+-+++++-+---+-++++-++++++
+-----+-+-+-+-++-+--++-----+--++++--++-
+--++--+---+++-+-+-+++-++------+-++++++
+----+--+--++----+++----++---++-++-++-+
+----++--+---+-+++++---+-+-++-++-++-+++
+----+-++-+--++++-+--+--+++----++-+++++
+++-+-----+-+---++++-++-+---+-+++-+-+++
\end{verbatim}
\end{minipage}
\end{equation*}

\vskip 0.5in

\begin{eqnarray*}
\text{Determinant:}&& 2^{38}\times9^{17}\times1241\\
\text{Theoretical bound:}&& 2^{38}\times9^{17}\times80\sqrt{357}\\
\text{Fraction of bound:}&& 0.821009
\end{eqnarray*}
\vfill

\begin{center}
Welsh~\cite{we} earlier reported a record of
$\,2^{28}\times9^{17}\times1197\,$ (79.2\% of bound).
\end{center}

\vfill
\eject
\phantom{.}
\vfill
$$\mathbf{n=45:}$$
\medskip

\begin{equation*}
\begin{minipage}[t]{3.6in}
\baselineskip 0pt
\begin{verbatim}
--+++--++++-++++++-+++++++-+-+++-+++--++-++++
-++-+++-++-++--+++++-+++--+++-+++++++++++++--
-+-+-++++++++++++-+-++-++++-++--+--++++-+++++
-+++++++--++-++--++++-+-+++++++++++-++-++--++
+-+-+-+-+-+-+++--+++++-++++-+----+++++-+++++-
++--++++-+--+---+-+-+-++++++-+-+++++-+-+-++++
+-+---+++--+-++++-+-+++++-++-+++--+-+++++-+-+
++-+---+++-----+-+++++-+++-+++++++-++++++--+-
+-++-++++++++---++--++---++---+++++-+-++++-++
+--++++++-+-++++-+----++-++++++++-+++++--+---
+--+++----++-+-++-+++++--+--+++-++++-++++++-+
+++--+---+++-+++++-+---+-+++++-+-+--++++-++++
++++----++-++++---++--+-+++++-+++-++--+-+++++
+---+-++-++++-++++-+-++-+-+++++-++-++---+-+++
+-+-++-+++-+-++++++---++++----+-++---+--++-+-
++++--++++--++-+++++-++--++-++---++--+------+
+++++-++-+++--+++----+--++-++--+--++-+-+++---
+++++-+-+-+---+++--++-++--+----+++-+-++-+--++
+-++++-+++++----+-++--+++-+-++----+++-++---+-
++---++++--++++-++-++-+-+---++-+-+-+--++++---
++-+-++--++--++-+++--++-+-----++-++++++---++-
+-++-+++-+++++-+--+++-++-+-+---+-+-++---+-+--
++-+++-+---+++++-++--+-+--++-----+++--+++--++
+++++-+----++-+-+-+-++++-+-++++--+--+-+--+-+-
++-++++-+--+---+-+-+-+++++---+-+--+-+---+++++
+--+++--++-++++-+--+++-++--++--++++-++------+
++--+++-++++++-+---++++-++++--+------+++---+-
+-++-++--++-+-+--+++-++++--+-+--+----+++++--+
+++--+----+++-++++++++-++++--++++-++---------
+++-+-+-++++-+-+-++-+---+--+-+--+++++-+--+--+
++++++--++---++-++--+++--+++-+--+--++--++-+--
++++-+-++-++-----+--+++++-+++-+--+-+-+---++-+
++-++----++-++-++++-+-+++-+-+-++----+--+++--+
++-+---++-+++-+++++++-+----+----+-+-++-+-+++-
++-++-++++++-++---++-+-+--+--+++++-----+-++--
+++-++-++-+-+-++--+--++-++--+--+++--+-++--+-+
++++++-+-+--+-++---++---+-+--++--++-+++-+++--
+++++-++---+++--++-+---+++----+-+--+++++--+-+
++--+++++++---+-+++++--+-+-++-+---+---+-+-+-+
+++--+++-+---+++---+++++----+-+-+-+++--+-+-++
+-++++++-----+-+++++++--+-+++--++-----+--+++-
+++-++-++-+-++--+-++-+-----+-+++---+++--++-++
+-+++++-++-++-++-++-+-------++++---+-+-+--+++
+++-+--+-+++++---+--++++----++-++-+--++-+-++-
++++-++-+-+-++-++------++--++++-+++----++-++-
\end{verbatim}
\end{minipage}
\end{equation*}

\vskip 0.5in

\begin{eqnarray*}
\text{Determinant:}&& 2^{44}\times 11^{21}\times 83\\
\text{Theoretical bound:}&&2^{44}\times 11^{22}\times\sqrt{89}\\
\text{Fraction of bound:}&&0.799817
\end{eqnarray*}
\vfill

Farmakis \& Kounias~\cite{fk} published an earlier record of
$\,2^{44}\times11^{21}\times81\,$ (78\% of bound).

\vfill
\eject

\phantom{.}
\vfill
$$\mathbf{n=47:}$$
\medskip

\begin{equation*}
\begin{minipage}[t]{3.6in}
\baselineskip 0pt
\begin{verbatim}
+++--++--++--++--++--++--++--++--++--++--++--++
+++-+-+-+-+-+-+-+-+-+-+-+-+-+-+-+-+-+-+-+-+-+-+
++++--++--++--++--++--++--++--++--++--++--++--+
--+----+++----++++-+++-+++----+---+---++++----+
-+-----++-+--+-++-+++-+++-+--+---+---+-++-+--+-
+------+-++-+--+-+++-+++-++-+---+---+--+-++-+--
+++-----++++----+++-+++-++++---+---+----++++---
--+---+----+++----++++-+++-+++----+---+---++++-
-+---+-----++-+--+-++-+++-+++-+--+---+---+-++-+
+---+------+-++-+--+-+++-+++-++-+---+---+--+-++
++++--------++++----+++-+++-++++---+---+----+++
--++++----+----+++----++++-+++-+++----+---+---+
-+-++-+--+-----++-+--+-++-+++-+++-+--+---+---+-
+--+-++-+------+-++-+--+-+++-+++-++-+---+---+--
+++-++++--------++++----+++-+++-++++---+---+---
--+---++++----+----+++----++++-+++-+++----+---+
-+---+-++-+--+-----++-+--+-++-+++-+++-+--+---+-
+---+--+-++-+------+-++-+--+-+++-+++-++-+---+--
++++----++++--------++++----+++-+++-++++---+---
--+---+---++++----+----+++----++++-+++-+++----+
-+---+---+-++-+--+-----++-+--+-++-+++-+++-+--+-
+---+---+--+-++-+------+-++-+--+-+++-+++-++-+--
++++---+----++++--------++++----+++-+++-++++---
--+---+---+---++++----+----+++----++++-+++-+++-
-+---+---+---+-++-+--+-----++-+--+-++-+++-+++-+
+---+---+---+--+-++-+------+-++-+--+-+++-+++-++
++++---+---+----++++--------++++----+++-+++-+++
--++++----+---+---++++----+----+++----++++-+++-
-+-++-+--+---+---+-++-+--+-----++-+--+-++-+++-+
+--+-++-+---+---+--+-++-+------+-++-+--+-+++-++
+++-++++---+---+----++++--------++++----+++-+++
--++++-+++----+---+---++++----+----+++----++++-
-+-++-+++-+--+---+---+-++-+--+-----++-+--+-++-+
+--+-+++-++-+---+---+--+-++-+------+-++-+--+-++
+++-+++-++++---+---+----++++--------++++----+++
--++++-+++-+++----+---+---++++----+----+++----+
-+-++-+++-+++-+--+---+---+-++-+--+-----++-+--+-
+--+-+++-+++-++-+---+---+--+-++-+------+-++-+--
+++-+++-+++-++++---+---+----++++--------++++---
--+---++++-+++-+++----+---+---++++----+----+++-
-+---+-++-+++-+++-+--+---+---+-++-+--+-----++-+
+---+--+-+++-+++-++-+---+---+--+-++-+------+-++
++++----+++-+++-++++---+---+----++++--------+++
--++++----++++-+++-+++----+---+---++++----+----
-+-++-+--+-++-+++-+++-+--+---+---+-++-+--+-----
+--+-++-+--+-+++-+++-++-+---+---+--+-++-+------
+++-++++----+++-+++-++++---+---+----++++-------
\end{verbatim}
\end{minipage}
\end{equation*}

\vskip 0.5in

\begin{eqnarray*}
\text{Determinant:}&& 2^{46}\times 11^{18}\times 15^4\times60\\
\text{Theoretical bound:}&&2^{46}\times 11^{20}\times19^2\sqrt{5665}\\
\text{Fraction of bound:}&&0.923897
\end{eqnarray*}
\vfill

Our conference-matrix construction produces this example (see
Sec.~\ref{sec:methods}). Welsh~\cite{we} earlier reported a
record of $\,2^{46}\times11^{21}\times1896\,$ (76.8\% of bound).

\vfill
\eject

\phantom{.}
\vfill
$$\mathbf{n=53:}$$

\begin{equation*}
\begin{minipage}[t]{4.2in}
\baselineskip 0pt
\begin{verbatim}
++------+-+----++--++---------+---+--+---------+-+-+-
+-+-+------+-----+---+-+-------+----+--++--+----+-+--
+-----++----++----------+--+++--+--+------+-+-+------
+--+-+---+----+---+---+--++------+----+--+---+------+
-+---++-+-++------+-++-+-+--+-+-+-----++----+-+---+-+
-+++---+-++--+--+++--++--+-+----+-+-+------++--+-----
--+--++-+++----+-+++---++----------+++--+++-++-------
---+--+++-+-++-+--+-------++---+-------+++-----++-+++
--++++-+----+-----+-+-+++---+--+--+--++---+----+-++--
---+------+++--++--+-+++--++++------+-+-+----++--+---
-+---+++--+--++--+-+--++--+--++--+--------++----+++--
---++-+----+-+++++-++-+-++--+---+----+-+-+------+----
--+++---+-+-+-+--+--+---+-+--++-+--++-+----+-------++
-+---+-+-+--+---++--++-----+-+++---+-+++-+---+--+----
-+--+--+-+-+-+++--+----+-+---+-+---+++--------+--+-++
--+-++--++----+-+--+-+---++++------+------+---+++-++-
-++--++-+--++---+-----+-++---+-+-++-----++-+--+----+-
--+-+-++-+-+---+---++-----+++-++-++--------+++------+
-+-++---++-++++------+--+-----+--+-----++-+-++-+-+---
---+--+--+--++-+++---+-+-+--+-+--+++-++-+-++------+++
--++--+-+----+---+++-----+--++++-++++-++-----+++++---
-+-++--++------+-+-+-+++--+----+++++--++-++-+-+----+-
-+----+-----+-++--++++-+++-+---+++--+-+---++-+-++--+-
---+-+-+--++---++-------+++---++++-+++-+--++-+++--+--
-+-++-+--+------+-+-+--++-++--+-+++++---++----+-+++--
-+--+-+--+++------+---+---+-++--+-++-++++-++-+-++--+-
-+--++----+-+--+-+----+-++-++----++++--+----++--+++++
-++------+--+--+-+++-+--+-+-++--++---+-+-+-+--++-++-+
--+-+-+---+---++-----++--+-+-+++++-+-++-++--+--+-++--
-+++-+----+--+++----++------+--++-++----++++-++-++--+
-+-+----+-----+-+-+++-+----+++-+-+-+++-++-+++-----+-+
-+--+++-----+++-++-+------+----++-+-+++-+---++++--+-+
-----+---+--++-----++-++-++-+-+++--++--+++-++--+-+-+-
----+---+---++-+--+-+++--++--++---+-++---++++++-+-+--
------+++-----+--+--++++++++------++-+-++--+-+++-+--+
---++++++-------+----+-+----++--++--++---+-+++-+++-++
---+-++--+-+--++-+--+------+-+----+-+-++-++++-++-+++-
-++----++--+-+---+--+----++++---++--+++-+++--+---+++-
----+---+++--+--+--+---+++-+-+-++-+---++-+++-+---++-+
--------+++++-+--+----++---++-+++++-++---++---+++---+
--++---+-++-+-+----++--+-+---+--+++--+-++---+++-+-++-
-+++-------+-+-----+--+++--+--+----+-++--+-++++++-+++
----++-+--+--++---++-+--+--++-+--++-++++++-+--+----+-
--+--+--+--++++++-+----+--++----++++-+++---++---++---
----+--+--+++---+++++----+-------+-+--+-+++++-++++--+
-------+++----+++++---+-+---+-+++-----+-+--++++-++++-
-----+-+++-+-+-+---++++-+----+--+++++-+-+------++-+-+
--+---++---++-+-+-++-++-------+-+-+++--+-++--+---++++
----------++-++-+++-++-++-+-++-+-+++-----+--++-+--++-
--+---+--++--+--+---+++-+-+----+-+--++++--+-+-+-++-++
---+-+-----+-----+++-+--++++-++++-+--+--+-+-+---++-++
-++----+------+++------++++-+++---+-+-+++++-+--++---+
--+-++-------+-++++-+-++---+-++-++-----++-+--+++---++
\end{verbatim}
\end{minipage}
\end{equation*}

\vskip 0.5in

\begin{eqnarray*}
\text{Determinant:}&& 2^{52}\times 13^{25}\times 105\\
\text{Theoretical bound:}&&2^{52}\times 13^{26}\times\sqrt{105}\\
\text{Fraction of bound:}&&0.788227
\end{eqnarray*}
\vfill

Farmakis \& Kounias~\cite{fk} published an earlier record of
$\,2^{52}\times13^{25}\times104\,$ (78.1\% of bound).

\vfill
\eject

\phantom{.}
\vfill
$$\mathbf{n=63:}$$

\begin{equation*}
\begin{minipage}[t]{4.5in}
\small
\baselineskip 0pt
\begin{verbatim}
+++--++--++--++--++--++--++--++--++--++--++--++--++--++--++--++
+++-+-+-+-+-+-+-+-+-+-+-+-+-+-+-+-+-+-+-+-+-+-+-+-+-+-+-+-+-+-+
++++--++--++--++--++--++--++--++--++--++--++--++--++--++--++--+
--+-------+---++++-+++----++++-+++----++++-+++----++++----+---+
-+-------+---+-++-+++-+--+-++-+++-+--+-++-+++-+--+-++-+--+---+-
+-------+---+--+-+++-++-+--+-+++-++-+--+-+++-++-+--+-++-+---+--
+++----+---+----+++-++++----+++-++++----+++-++++----++++---+---
--++++--------+---++++----+---++++-+++----+---+---++++-+++-+++-
-+-++-+------+---+-++-+--+---+-++-+++-+--+---+---+-++-+++-+++-+
+--+-++-----+---+--+-++-+---+--+-+++-++-+---+---+--+-+++-+++-++
+++-+++----+---+----++++---+----+++-++++---+---+----+++-+++-+++
--++++-+++-----+++----++++-+++----+---+---+---+---++++-+++----+
-+-++-+++-+----++-+--+-++-+++-+--+---+---+---+---+-++-+++-+--+-
+--+-+++-++----+-++-+--+-+++-++-+---+---+---+---+--+-+++-++-+--
+++-+++-+++-----++++----+++-++++---+---+---+---+----+++-++++---
--+---++++----+----+++----++++----++++-+++----++++----++++----+
-+---+-++-+--+-----++-+--+-++-+--+-++-+++-+--+-++-+--+-++-+--+-
+---+--+-++-+------+-++-+--+-++-+--+-+++-++-+--+-++-+--+-++-+--
++++----++++--------++++----++++----+++-++++----++++----++++---
--+---+---++++----+----+++----+---++++-+++-+++----++++-+++----+
-+---+---+-++-+--+-----++-+--+---+-++-+++-+++-+--+-++-+++-+--+-
+---+---+--+-++-+------+-++-+---+--+-+++-+++-++-+--+-+++-++-+--
++++---+----++++--------++++---+----+++-+++-++++----+++-++++---
--++++-+++----++++----+-------+---++++-+++-+++----+---+---++++-
-+-++-+++-+--+-++-+--+-------+---+-++-+++-+++-+--+---+---+-++-+
+--+-+++-++-+--+-++-+-------+---+--+-+++-+++-++-+---+---+--+-++
+++-+++-++++----++++-------+---+----+++-+++-++++---+---+----+++
--+---++++----+---++++-+++--------+---+---++++-+++-+++----++++-
-+---+-++-+--+---+-++-+++-+------+---+---+-++-+++-+++-+--+-++-+
+---+--+-++-+---+--+-+++-++-----+---+---+--+-+++-+++-++-+--+-++
++++----++++---+----+++-+++----+---+---+----+++-+++-++++----+++
--+---+---++++-+++-+++-+++-+++-----+++----+---+---+---+---++++-
-+---+---+-++-+++-+++-+++-+++-+----++-+--+---+---+---+---+-++-+
+---+---+--+-+++-+++-+++-+++-++----+-++-+---+---+---+---+--+-++
++++---+----+++-+++-+++-+++-+++-----++++---+---+---+---+----+++
--++++----++++----+---+---++++----+----+++----++++-+++----++++-
-+-++-+--+-++-+--+---+---+-++-+--+-----++-+--+-++-+++-+--+-++-+
+--+-++-+--+-++-+---+---+--+-++-+------+-++-+--+-+++-++-+--+-++
+++-++++----++++---+---+----++++--------++++----+++-++++----+++
--+---++++-+++----+---+---++++-+++----+----+++----+---++++-+++-
-+---+-++-+++-+--+---+---+-++-+++-+--+-----++-+--+---+-++-+++-+
+---+--+-+++-++-+---+---+--+-+++-++-+------+-++-+---+--+-+++-++
++++----+++-++++---+---+----+++-++++--------++++---+----+++-+++
--+---++++-+++-+++----+---+---++++-+++----+----+++-+++----+---+
-+---+-++-+++-+++-+--+---+---+-++-+++-+--+-----++-+++-+--+---+-
+---+--+-+++-+++-++-+---+---+--+-+++-++-+------+-+++-++-+---+--
++++----+++-+++-++++---+---+----+++-++++--------+++-++++---+---
--++++-+++-+++----++++-+++----++++----++++----+-------+---+---+
-+-++-+++-+++-+--+-++-+++-+--+-++-+--+-++-+--+-------+---+---+-
+--+-+++-+++-++-+--+-+++-++-+--+-++-+--+-++-+-------+---+---+--
+++-+++-+++-++++----+++-++++----++++----++++-------+---+---+---
--+---+---+---++++----++++----++++----++++----++++-----+++-+++-
-+---+---+---+-++-+--+-++-+--+-++-+--+-++-+--+-++-+----++-+++-+
+---+---+---+--+-++-+--+-++-+--+-++-+--+-++-+--+-++----+-+++-++
++++---+---+----++++----++++----++++----++++----+++-----+++-+++
--++++----+---+---+---++++-+++-+++-+++----++++-+++----+-------+
-+-++-+--+---+---+---+-++-+++-+++-+++-+--+-++-+++-+--+-------+-
+--+-++-+---+---+---+--+-+++-+++-+++-++-+--+-+++-++-+-------+--
+++-++++---+---+---+----+++-+++-+++-++++----+++-++++-------+---
--++++----++++-+++-+++----+---+---+---+---++++-+++----++++-----
-+-++-+--+-++-+++-+++-+--+---+---+---+---+-++-+++-+--+-++-+----
+--+-++-+--+-+++-+++-++-+---+---+---+---+--+-+++-++-+--+-++----
+++-++++----+++-+++-++++---+---+---+---+----+++-++++----+++----
\end{verbatim}
\end{minipage}
\end{equation*}


\begin{eqnarray*}
\text{Determinant:}&& 2^{64}\times15^{24}\times 19^7\\
\text{Theoretical bound:}&&2^{64}\times15^{28}\times12^3\times2\sqrt{33}\\
\text{Fraction of bound:}&&0.889364
\end{eqnarray*}
\vfill

Our conference-matrix construction produces this example (see
Sec.~\ref{sec:methods}). We find no prior records in the
literature.

\vfill
\eject

$$\mathbf{n=69:}$$

\begin{equation*}
\begin{minipage}[t]{5in}
\small
\baselineskip 0pt
\begin{verbatim}
+-+++--+++-+-+++++--+++-+++++++++++--+++-+++++++++++++-++-+-++-+++-+-
+++-++++++++++----+++++++---+-+--++++++++-+++-++++++-++++++++++---+++
+++++++-+++-+-++++++-+-+-+++-+-+++-+++++++-+-+-+----++++-++++++++++++
++-+-+++--+++++++++++-+++++++++++-+++---+++-+++-+++++-+-++-+--+++++-+
-++---+++-+++--+--+++-++--++--+++++--+---+--+++-+--+-+++--+-++++++-++
--++--+++++--++--+++--+++-+-++++-+----++++--++-+-++----++-+++-++-++-+
-+-+-++--+-++-+-++-+++--+---+-+-++--++++++--+++--+-++-+-+-+++++++--+-
---+++++-++-+--++++--+-+--+++-+---+----+--+-+++++++-++++-++++--++-++-
-++-++--++-++++--+-----+-++-++++--+--++-+---+++++-+-+++-++-+-++-++-++
-----+-++++++---+--+--+-++++++-++++-+-++---+-+++++-++-++-+-++++--++--
-+--+-+-+-++++-+++-++++---++++--++--+-+-+-++++-++-+--+--+--+++-++-+-+
--++-+------++++---+++-+-++++++-++--+-+--++++-++-+++-++--++-+-+--++++
---++++++-----++---++++--++-+++-++++-+-++------++-+++-+++++--++++-+-+
-++++-----+-+-++--+++---+--++-++--+++-++++-+++-++-+++--+++--++--++++-
--+-+++--+++--+-+++-+--+---+++-+++-+-++-+-+++-+-++++---+-++--++++++--
--+-+--+-++-+-+-+-++++--+-+--+-+++++++--+++---+++++-+++-+-+++---++--+
-+-+-+--++++++-++--++-+--+---+-+--++-++++++-+--+-+++-+++--++---+++++-
-+--+-+--++-++++++--++++++++-+----++---++--+-++--+++--+++-+-+++--+-++
-++++--+-+-++---++++--+--++++-+-+--+---+++++---+++-+-+--++++-+++-+-++
-+--++++---+-+-++-+--+--++-++-++-+++-++-++-+-+++--++++----+++-++-++-+
-+++--+++-----+-++-++++++--+-+----+--++-++++-+++++-+++-+-+-+--+-+-+++
---++---++++-+-+-+++-+++-----++++++++----+-+-++--+++++--++++-++-+-++-
---+-++++-+++--+-+--++-+++---++++--+--+++++++-+-+---++-+++++++---+--+
--+-+-++--+++++-++----++-+----+-+++-++++-+++------+++-++++-++--++++++
-+++++-+++-+-+----+--++-+-++-+--+--++-+-+-+-+++----++-+++++-+--++++++
-+-----+++---+++++++-+-+-++-+--+-++++-+++++++-+-+--+-+++++--+-+++----
---+-++---+++++---+--++-++++---+-+-+++++-++--++++++--+-+++---++++--++
-+-+--+-+-+--++-+-+---++++--+++-++++++--+-++++-+++--+++--++-++-+-+-+-
--+--+---+-+-+++-++++-++++++--+-+-++++-++--++-++-+--++-+---+++-++-+-+
--+-+-++++-+---++--++--+++-+++-+---+++-+-++-++---++-+++-++--++++--+++
-+---+-+++--++++-++-+++++--++--++---++---+++---++-+-+-+++-++-+++-+++-
--++++-++-+--+-++++-++--+++---++-+--++-++-++++--++-+--+----+-++-+++++
-++++-+--++--+-+-+--+---++-+---++++-+-+-----+--+++-+++++++++--++----+
-++-++--+---++-+++++---++----++-++-+--++-++--++++++++-++-----+-+----+
-+-++-+-++++-----+++-+--++-++++-+-+-++-+-++++-++--+---++------+-++--+
-++++++-+++++++----+-+-++++-+--+--+------+++----++-++-----+-++-++-+-+
--+++++++--++++++--+---++--+-+---++++--+---++-+++-------+-++--++++-+-
-+-+++---+-+--++--++++++---+-+++-++--+++--++-+--++----+-++-++--+-+--+
-++--+++-++----+---+-+-++++-+++-+--++++----+-+--++++-+-++--+---+++-+-
-++++----+++++-++---+++++-+-+++--+--++-+++----+-+---++-+-+----++-++--
-++++--+--++--+-+--+++++-++---++---++-+---++++++--+-+-++--++-+++-----
-++--++----+--++++++-++-+-+-++-+-++-+--+--+-+-----++++-+-+++-+---+-++
--+-+-+-----++---++++++-+++-++-++-++-+++-+-++++-+-----+--+++---+--++-
-+-+++-+--++-++-+++++--+-+--++-++---++-+-----++++--+-+-++-+-+-----+++
-++-+++++-+--++--+--+++----+++++---++--+-+----++-+-+-++--+-+++-+++---
--++++--+-++--++-+----+++-+++---+-+++++++++--+--+-++-++---+++-+------
----+++-++++--+-++++++++++++--+--+-+--+--+--+--++--++++-+--------+++-
--++--++-+-+++-+++++-+++++-+-+++----+++-+------++-++--+--++-++--+----
-+++-++-++---+--++-++----+++--+--+++++---++--++-+-+----+++++++---++--
-++-----+++-+-+++----+++++-++-++++-+-+-+--+-++++---+----++++----+-+-+
-+++++-++-+-+-++++++-++-++------+-+--++-----+-+--++--+--++--++++-++--
-++++++-+--+---++-+-+-+++++-+--+++------++-+--++-++---++++-++---+--+-
-++--+-++-+++----++-++-+-+++-+---++-++-+++-+++-+-++++---+++----+-----
---++-++++-++++++-+-++-+--+-+---+-+++++--+--++-+-+-+---+-+-+-+---+--+
--++++++-++-++++--++--+--+-+++---+---+--+++++++----++++++--+-+-------
-+++--++--++-+-+-+-+-+-+--+++--+++++-++++-+---++-+--+-+-+----+---+++-
-+--+--++-++++++-++++---+++++++--+-+-++---+------+--+--+-++++-+-+--+-
-++--+++-+++-++++-+-+---+-+--++-+-+---++-+++-+-+-----++-+++--++-+-+--
-+-+++-+----+---++----+++-++++++-++++++--+-++----+---++++-+--+--+-+-+
-+--+++++++--++++--+--+---++--+++---++++++-++-+-+++-+----++--------++
-+++-+-+-++--+--+-+-++++-+-+++--++++--++-+--+---+-+-+------++++-+--++
-+-+++++-+--++---+-++-+++-+----+++-+-+-+--++++-+--++++---+--+-+-++---
--+-+++-+---++--+-++++++-+-++-+++-+-+-+++-+--+---+--++-++-+---+-++---
-++-++++-+-++-++-+---++--+---++++++++---+++-++-+++-----+----+-+---++-
--++-+-++++++++-++-+++--+--++-++++++----+--+-+----+--+++-+-----+---++
-+++--++++-++-++-++-+-+--+--+----+-++-++---+-++-+++-+++---+-----+++-+
--++-++-+++-+--++++-+-+---+-++++--+++++---++--++---++---+---+-++---++
-+++--+-+--++++-+-+--+----++-++++--+-+-+---++---+++++++++--++-+---+--
-+-++-++++--+-+---+-+--+++++-++++++-+-+++-+---+---++-+-----+-+-+--++-
\end{verbatim}
\end{minipage}
\end{equation*}

\begin{eqnarray*}
\text{Determinant:}&& 2^{68}\times 17^{33}\times 155\\
\text{Theoretical bound:}&&2^{68}\times17^{34}\times\sqrt{137}\\
\text{Fraction of bound:}&& 0.778973
\end{eqnarray*}
\vfill

Farmakis \& Kounias~\cite{fk} published an earlier record of
$\,2^{68}\times17^{33}\times153\,$ (76.9\% of bound).

\vfill
\eject

$$\mathbf{n=73:}$$
\begin{equation*}
\begin{minipage}[t]{4.75in}
\footnotesize
\baselineskip 0pt
\begin{verbatim}
+------------------------------------------------------------------------
-+++---+-+-+---+-+-+-+-++-++--+--++-+---+--------++-++----+-++++-+-------
---+-+++---++++++++++----+-----+--+++-+---+--+--+---+--+------+++------+-
-+++++----+-++--------+----++-++-++-+++-+-++----++-+--+--+-+---+------+--
--------++-++-+-----+++-++----++-+-+++++-+-+---++++++--+--+--+------+---+
-------++++++---+-+-+-+-+-+-+---+-+---+++--++---+-++-+-----+--+++++--+---
--+++-+---++---+-+--++++----++-++--++--++--+---++-----+-----++--++--++-+-
-+++-++------++-+---++-+-++-----++----+--+-+----+-++-++++-+++---++-+-----
-+++---+++-+++-------+-+-++-++--+---++----+++-++----+--+---+---+--+++--+-
---+++--++-+-+-+--++++----++--+---+----+-++++-+++-+---+-+++----------+-+-
-+---++---+-+--+---+-++-+-+-++--+-+++-+----++++--+--+++++-+----------+--+
--+----+++--++-+++-----++--+----++-+++-+++---+--+----+++++-+---------+-++
-+---+++------+--+++-++--++-+-+---+-++-+++-----++----+--+--+-----++-+-+++
--+-++-+--++-++-++--+---+-+--+-+----+----+++++-+-+---+---+++-+---++-----+
-+-----++--+-+-+++++++------++-+++---++----+-----++---+-++--+-+---+-+-+-+
-++---+++--++-+--+--+-+-++-+--+-+-+------++-+-+--+----++-+--+---++--+++--
----+-+-+-+++-+-+--+---+----+-+-++--+---+++--++++++-+-+-+-----+--+++-----
---++--+--+--+-+--++++--+--++-+++----+-+-+---++--+-+++-+---+----++-++----
-+---++-++----+++-+-+--+---++--++--++---+-+-+-+--+++-+----++-+------+-++-
--+---+-+-++++----++-+--+++-+-++-+-+----+-+-+------+-+--++--++--+--+---++
--++-+----+-+++-+++----+++--+++-+-+--+-+------++-+++------+-+-+------+-++
----+++++----++-----+---++++++---+--++-++--+-++---+-+---+---++-++----++--
---++----+--++++-+--------+--+--++++--+++++-+--+-+-+++--+---+---+-+-+-+--
-++-+----+--++--+-+++-++-+-++--+------+-++--++-+----+---+-+-+--+-+--++--+
---+-+-++-+------++-+-++++-+-+++++-+--------+--++---+++-+-+----+--++--+--
-+++----+--+--++---------+-+++-+--+--+++++--+++-+----++---+--++-+-++----+
-+-++-++-+--+---+--++--++-+-+++----+-++--+------++-------++--+--+-++-+++-
-+---+-+--++--+++-+--+++-----+-------+-++++-+---++-++--+++--++-+---+-+---
-+--++-+-+--+---+++--+---+---++--+-++--++--++++----+--+--++---+-++-++----
-++------++---+-++--+++---+++--+-+++---+-----+++-++----+++-----++-++---+-
--+--+--+++------+++---++-+--+-+--++-++-++++-++-+-+----+----+----+-++-+--
--+-++-+----+--+---+-+++++-----++-----+++-+---++--+----+++++-++-+-+---+--
-+---++--+-+-+----+---++------++++++-+---+-+-+++----++---+-++++-+-+--+---
--+-+++--+++---++--++--+-++---+-+++--+-+----++-----+--++-----+-+--+-+-+-+
-++---+--+++-+++--++----++-+-+-----++-++-------+--+++-+--+-+-----+++-++--
---+-++-+----+--++---++++-+++------+--++-++---+----++-+--+---+++-++-+---+
-+-++-+++++--++------+++-----+++-+----++--+-++----+--+-+------+--+---++++
-+--+-+-+-------++-+----+--+-++++-+-+-++-+-++--+---+---++--+++++---+---+-
--+-+-++++----+---++-++---++-+--+--+-+----+--+-+++-+--+---+++-+++-------+
-+--++--+--++--+++--------+-+-++---+-+-+--+----+--+-++++--+++--+++-+-++-+
-++++---+-------+-+-++--++--------++-+--+-++-+-+-+-++++++----++--++--+++-
---+--++-++++++-++-+--+---++---+-----+--+--+--++---+-++++-+--++----++++--
--++----+-----+++-++--+---+--+++++--+---++-+--+----++--+-++----++++--++++
---++-++------++--+----++++-+--+-+++----++-++--+-++---++-+-+--++---+++--+
--+++-+--+-+----+-+-+-+-+---++---++-+--+--+---+-+---++-++++++-+----++-+-+
---+----+++---+--+++---+-+--+-+-----+-++--++-+---+---+++++++++-++-+-++---
--+--+-+----+--+--+-+-+---+++++--+--+-+--++--+-+--++-+++----+++---++++-+-
---------+--+-++-+--++-+-++++++++-+-------++-+--+--++---++-+-++--+-+-++-+
-+-+--+++-+-++-+----+--------++---++----+----++++-++---+++++++-+-++-+---+
-++-+--------+++-+-++-+++-----+----+++----+++-+-+-++-+--+-----+++-++++--+
--+-+--+--+++-+-+-+--+-+---+-++---++--+--+-+--+---+-+++-+--+-+-++---+-+++
----+-++-++------+--++---+-----++-+-+++-+++---+---+++++--++-+--+--++-+-++
-+--++--+-++-----+-++--+-+++-+---+----+-+-----++++-+++-+-+-+--+-+---++-++
----++--++--+++++--+-++-++---+-+-++-----------+-++---++----+++-+-+++++-+-
-+-+-+----+++-+----+++--+--+----++-+---+++++-+------------+++-++--+++++++
------+-----++---++++-++---+-+---++--+-++-+++----++-+-++--++-+--++++---++
----++-+++++-+-+------++-+-++-----+++----+-+-----+-+---++-+-+-+-++++--+++
--+---++--+----++--------+++--++-+-+-+++--+++-+-++--+---+-+++-+--++-++-+-
---+-++--+-+---++-----+-++-+-++-+----++-+-++-+-+-++--+--++-----+++-++--++
--++++-++--+--+--+-+--+-----+---++++-++-----++----+++---++++-+---+-+++-+-
----------+--+--+--+-+---+-+++--++++++---++-+--++-+--+-+-++--+++++--+++--
---+-------++------++-++-++--+-+---+++-+-+--+++--+++-++-++-++-++-+----++-
-+-+++---++-+-+---+-----+-++---++---++------+---+-+-+-++++--+++-+++-+--++
-++-+-+-+-+-++-+-++---+--++-----+------+-+-+-+--+++-++---++--+++---++-++-
--++-++-+++-+--+-++-++-----------+--+++--+--+-++-+------+--+-+++++-+-++-+
---++------+-+--+++--++--+++--+-+--++-+-+---+++-+++----+---+++----++-++-+
----+-+----++-++-++--+-++--++----+---++----++++++--+-+-+-++-+--+-++---++-
-+-+---+-++----+++-+--+-++--+----+-------++++++++-+++-+----+++--+-----+++
--------+-+--++++---+++++-+---+---+-+++-+---++-+----+-+--++++-+-+--++-++-
--++-++++---+---+--+-+-+---+---+--+-+--+---+++-+++++++---+--+---+-+++++-+
--+--+-+-+---++----++-------+-+-+--+--+++--++-++++--++++-+--++++-+-+--+++
-+-----+---+-+----+----++-+----++++-+-++--+--+++++-+--+-+-+--+-+++-++++++
--+--++--+-+-+---+---+--+---++++--------++--++--+++++-+++-++--+++-+-++++-
\end{verbatim}
\end{minipage}
\end{equation*}

\begin{eqnarray*}
\text{Determinant:}&& 2^{72}\times 18^{35}\times 163\\
\text{Theoretical bound:}&&2^{72}\times18^{36}\times\sqrt{145}\\
\text{Fraction of bound:}&& 0.752023
\end{eqnarray*}
Farmakis \& Kounias~\cite{fk}, published an earlier record of
$\,2^{72}\times18^{35}\times162\,$ (74.7\% of bound).
The lower right $\,72\times 72\,$ block above gives a Hadamard
matrix with excess 580. The largest excess previously known for
Hadamards of this size was 576~\cite{fk}.

\vfill
\eject

$$\mathbf{n=77:}$$

\begin{equation*}
\begin{minipage}[t]{4.5in}
\scriptsize
\baselineskip 0pt
\begin{verbatim}
+-++++++++-+-+-+++--+++++-++++++++-+++++++++--+++--++-+-+++--+++++-+--+++++++
+++++-+-+++++-++--++++++-+-+++----+++++++-++++++-+++++-+++++++-++-+-++-+++-++
++---+-+++++++++++++-++-++++--+++++++++-+++-+++++++-++++-+-++-+--++++++-+-++-
++++++++--+-+++-+++++--++++-+++++++----+-+-+++--++++-++++-++++++++++++++-++-+
---+-+-+----+-++++++++++--++++-++--++---+-++++-++-+++++-++-+++-----+++++-++++
-+++--+---+++--+-+++-++-+--++++-+++--+---+-++++++++-+-++++-+--++++----+++++++
-++---+----+-++-++---+-++++++--+-++++++---+---+-++++-+++++++++++--++-+++++-++
--+++--+++++---++-++----+++-+-++++++-++---+-++-+--++------++-+-+-+-+++-++++++
--+--++++--+--+-+--++++--++-+++--++-++++++--+--+++++++++--+-----++--+++--+-++
-+--+++-+++--+--+---++++++--++-++-+-+-+-++-++++++-+---++++-+-+-+-+-+++--++---
--++++---+++-+---+----++-++++--+++--++-++-+++++-++++++++--+++---++--+--++-+--
---+-+-++++-+---++--++++--++-+++-++--+--+++-+-+--+-----++-+++++++-+-++--+++++
-++---+-+++-+-+-+++-++--+-++----+++---+++-+++-+++--++-+--++-+--+++++++-+--+-+
-++-++---+++++--++-+++---+--+-+---++-+-+++++-+++---++++----++++++--+-++--+++-
----+-+-+++-++++-+-++------++++++-+++----++----++++++-++--+-+++-+++++--+++-+-
--+++--+-++-+++-++++-+++-++--+----++--++++-+--+++++--+-+-++--++--+-+---++-+++
-++-++--+-+--+-++-+-+-+---+++-+-+--+-++-++-++++--+-+-+-++++-+-+--++++-++-+-++
-++++-----+-++++--+-+-++++++--++++--++-+++----+++-+-++-------+-++-+++++-++-++
-++-++-+++-++++-++-+---+++-++++-+-----+---+-++++++--++--+++--+--+-+-+-+-+++-+
-+-+-+++++---+-++-+++----+--+-+++++-+--++--+--++++++++--+++++-++--+-----+-+++
--+--+--++++--++-+-++-+-+-+--++++-+++-+++-++-++-+++--+--+--+-++++++--+++----+
--+-+-+++--+-+-+++----+++-++-++---+++--+-+-+++-+--++++-+-+-++-+++-++++--+-+-+
--++-++-+-+-+--++-+--+-+++-+-+++-+-+--++-+++-+++-++++++-+-+-+++--+--+++-+----
-+-+-+-++----+++-++++-++++-++------+--+++++-+-++--++--++--++--+++++--++++++--
-++-+--+-+--------++++-+++-+-++-++-++++++++--+--+-++--++++-++++-+++-+-----+++
-+-++--++--++-++-----++--++-+++++++---+++----++----+++++++--+++-++++-+-++++--
--+++++--+--+-+----++++-+++-+--++--+++-+-++-++-+-+--+--++++-+-++-+++-++-+-+-+
---++++-++-++-+---++---++-+---+--+-++++-++-++-++++-++++++--+-+++--+++--+-++--
---+--++-+++-+++----++-+++-+--++--++--++++--++--++--++--+++++--+++-++-++-+++-
-+--+-++-+-+++-+--+-+-++--+-++-+-+++-+++--+++++++++---+-+-+---+-+-+++++---++-
----+-+++-++-+++++++++++-+----++-+--++-+--++-+++-+----++-+--++-+-++-+-++-++-+
---+++-+-+-++++-+-+-+-+-+--+-+--+++-+++--+-+-+-+-+-+-++--+++++-++++--++++--+-
-+-++++++-++++--++++-+--+----+-+---+++-++----++-+--++--++-++----+++++++++--++
-+---++++++----+-+++---+-++--+--++--++++-+++++--+---++-+-++-+++++--+-+++++-+-
-+-+++++++++-----+---++-++-+++--++-+---++---+--+-++--++-----++++--+++-++---++
-+-++-++++++++----+-+++-+-+++-+-----+-+---++----+-+-++-+--+++++--+---++--++-+
--+-++-+-+-++-++++++++--+--++--+-+--+-++---++-+---+--+-++---+-++++-++---++-+-
---+++--+++--+++-++--+--+++-++-+--+--++--+++------++++++++-----++-+-+-+---+++
--++-++-+-++++---++++-+-++++++-+-++++-++-+--+++-+--+-----+--++----+------+++-
-----++---++++-++-++-++++---+-+-+-+-+++++++-+-----+--++-+-+-++--+-++---++-+-+
-+++-+--++-+-++++-+-++-+---+++-+-++--+-+-++-++-++---++-+---+-++--+++---+----+
-++--++++++-+-+--++---++---++-++-+++++--+----+-+--++-+-++++--+-+++-+--+---+--
-+++---+++--++-+---+-++-+------++++++--+-++++-+--+-+-+++-+++-+--+--++-+--+--+
--+++--++-+---+-+-++-+++-+-++++++-+++++----++-+-+---+-+---+++-+++-++--+------
-++-++++-++---+++-+-+-+-++++-++---+-+--+--+-+-+--+-++-+++--+-+-----+--+++-++-
--+++++-++--++-+-+++++-+-++++-+--++-+-+-+----+---+----+--+-+--+-+--++++++---+
-+++-+++++-++++-+-+---++-++-----+-+++----++--++---+-+-+-+---+-++++--+--+-+-++
-+--++++----+++----+------++++++-+++-++++++++-+---+-+----+-+++-+-+--+-+++---+
-++-++-+---++--+-++--+++++--+-++--+-+---+++++---++-++----++--+++-++-++-+---+-
-+++-++--+++-++-+--+++++------++++-+--+----++-----+++--+-+---++-+-+-+++++-++-
-++++-+++-+--++++--+-+----+++----+-++---++++++--++---++-+----+-++++--+--+-++-
-+-++----+++-++-++++--+----+-+++-+--+-+++++-++---++++-+-+++---++---+-+---+--+
------++-+++++-++-+++--+++++++-++--+-+--+---+-+--+-+++++-+---+++++---+------+
-++-+++-++--++++-++--+++-----++++--+-+++----+---++-+-++----++--+-+---+--+++++
-+--++--+++++-+++--+--+++++++-++-+-+---+-+-+---++-----+++++++---+----+-+---++
----+-+-++--++++++-+-+++++++----+++--+-++--+++----+++---+-++-++--++---+--+-+-
-+-+----+---++--++-++-++++-++++-++++++++--++-----+--++-++------+-+-+++-+--++-
-++----+++++-++---++-+-++-++++-++-+-++-+++-+---+-+-+----+---+-+-----+++++++--
-+-++-+++-+-+-+-++--++-++-+-+-+++---++++-++-++++-+++-+---+-+--+-+------+---+-
-+++-+++-+--++-++--+-++-++++-+++-----++-+-++-+-+++++---+----+--++-++---+-+---
--+--+-++-+++++---+-++-++-+-+++-++-+---++-++-+---++-+-++-+++---+--++----++-+-
--+++-++++-++---+++-+----+-+--+++-++-+-++++++---+-+--+++++--------+--+++++---
-+-++-+--+-+--+++++++++++-+-+++-++++----+++--++-+--+-+-+--+-+--+----+-+-+----
-+-+-+-+--++++++-+--++---+++--+-+-++++++---+++-++--+--+++-+--+++----+---+---+
-+++++-+-++--+-+-+-+++-++-+--++--+++++------+-++--+-+++-+++-+----++--+-+-+---
--++--++-+++++++-+-----+-+--+++-++--+-+-+++++-++---+---++--++----+++-++----++
--+-+++++-++++++---++++--+-+-+--++---+---++---+-+-++-+--++++--++---+++-+--+--
--++++++--+--++++++---+-+---+---++++-++++-+----+++--+--+++-++++-+---++---+---
-++-+---+--+++-++++-++---++--+-+++-++-+-++--+--+-++-+--++-++++-++-----++--+--
-+++---++++++--+++-++-++--+-+--------+++-+------+++++-+-++-+++-+-++++-+-+----
-++++++++--+--++-+-++--++--+---++-+--++-++-+-++--++-+-++--+-++-----+-+----+++
-+++--+-+--+-++--++++-+-+++--+++-----+--+-+++--+-----++-+++++++-++-++---+--+-
-+-+++--+-+++-+++---+--+++---+---+++++--+-+++---++++-----+--+-+-++-+--+-+++-+
-++++++---+++--++--+---+--++-+-++-+-+-+++-+----+-----+++-+++--++-++-+-+--+++-
-++---++-++-+-++++--+-+-++--++-+-+-++++-++-+-+------+++---++--+---+-+--++++-+
-+-+-++--+-+--+++++--+---++++++-+--+++-+-+-+--+--++------+++-+--++++++---++--
-+--++++--++-----+++++++-+++---+++++--+--+++---+-+-+++--+--+----++++----++--+
\end{verbatim}
\end{minipage}
\end{equation*}

\begin{eqnarray*}
\text{Determinant:}&& 2^{76}\times 19^{37}\times 177\\
\text{Theoretical bound:}&&2^{76}\times19^{38}\times\sqrt{153}\\
\text{Fraction of bound:}&& 0.753137
\end{eqnarray*}
\vfill

Farmakis \& Kounias~\cite{fk} published an earlier record of
$\,2^{76}\times19^{37}\times174\,$ (74.04\% of bound).

\vfill
\eject

$$\mathbf{n=79:}$$

\begin{equation*}
\begin{minipage}[t]{4.5in}
\scriptsize
\baselineskip 0pt
\begin{verbatim}
+++--++--++--++--++--++--++--++--++--++--++--++--++--++--++--++--++--++--++--++
+++-+-+-+-+-+-+-+-+-+-+-+-+-+-+-+-+-+-+-+-+-+-+-+-+-+-+-+-+-+-+-+-+-+-+-+-+-+-+
++++--++--++--++--++--++--++--++--++--++--++--++--++--++--++--++--++--++--++--+
--+----+++----+---++++-+++-+++-+++----++++----++++----+---+---+---++++-+++----+
-+-----++-+--+---+-++-+++-+++-+++-+--+-++-+--+-++-+--+---+---+---+-++-+++-+--+-
+------+-++-+---+--+-+++-+++-+++-++-+--+-++-+--+-++-+---+---+---+--+-+++-++-+--
+++-----++++---+----+++-+++-+++-++++----++++----++++---+---+---+----+++-++++---
--+---+----+++----+---++++-+++-+++-+++----++++----++++----+---+---+---++++-+++-
-+---+-----++-+--+---+-++-+++-+++-+++-+--+-++-+--+-++-+--+---+---+---+-++-+++-+
+---+------+-++-+---+--+-+++-+++-+++-++-+--+-++-+--+-++-+---+---+---+--+-+++-++
++++--------++++---+----+++-+++-+++-++++----++++----++++---+---+---+----+++-+++
--++++----+----+++----+---++++-+++-+++-+++----++++----++++----+---+---+---++++-
-+-++-+--+-----++-+--+---+-++-+++-+++-+++-+--+-++-+--+-++-+--+---+---+---+-++-+
+--+-++-+------+-++-+---+--+-+++-+++-+++-++-+--+-++-+--+-++-+---+---+---+--+-++
+++-++++--------++++---+----+++-+++-+++-++++----++++----++++---+---+---+----+++
--++++-+++----+----+++----+---++++-+++-+++-+++----++++----++++----+---+---+---+
-+-++-+++-+--+-----++-+--+---+-++-+++-+++-+++-+--+-++-+--+-++-+--+---+---+---+-
+--+-+++-++-+------+-++-+---+--+-+++-+++-+++-++-+--+-++-+--+-++-+---+---+---+--
+++-+++-++++--------++++---+----+++-+++-+++-++++----++++----++++---+---+---+---
--+---++++-+++----+----+++----+---++++-+++-+++-+++----++++----++++----+---+---+
-+---+-++-+++-+--+-----++-+--+---+-++-+++-+++-+++-+--+-++-+--+-++-+--+---+---+-
+---+--+-+++-++-+------+-++-+---+--+-+++-+++-+++-++-+--+-++-+--+-++-+---+---+--
++++----+++-++++--------++++---+----+++-+++-+++-++++----++++----++++---+---+---
--+---+---++++-+++----+----+++----+---++++-+++-+++-+++----++++----++++----+---+
-+---+---+-++-+++-+--+-----++-+--+---+-++-+++-+++-+++-+--+-++-+--+-++-+--+---+-
+---+---+--+-+++-++-+------+-++-+---+--+-+++-+++-+++-++-+--+-++-+--+-++-+---+--
++++---+----+++-++++--------++++---+----+++-+++-+++-++++----++++----++++---+---
--+---+---+---++++-+++----+----+++----+---++++-+++-+++-+++----++++----++++----+
-+---+---+---+-++-+++-+--+-----++-+--+---+-++-+++-+++-+++-+--+-++-+--+-++-+--+-
+---+---+---+--+-+++-++-+------+-++-+---+--+-+++-+++-+++-++-+--+-++-+--+-++-+--
++++---+---+----+++-++++--------++++---+----+++-+++-+++-++++----++++----++++---
--+---+---+---+---++++-+++----+----+++----+---++++-+++-+++-+++----++++----++++-
-+---+---+---+---+-++-+++-+--+-----++-+--+---+-++-+++-+++-+++-+--+-++-+--+-++-+
+---+---+---+---+--+-+++-++-+------+-++-+---+--+-+++-+++-+++-++-+--+-++-+--+-++
++++---+---+---+----+++-++++--------++++---+----+++-+++-+++-++++----++++----+++
--++++----+---+---+---++++-+++----+----+++----+---++++-+++-+++-+++----++++----+
-+-++-+--+---+---+---+-++-+++-+--+-----++-+--+---+-++-+++-+++-+++-+--+-++-+--+-
+--+-++-+---+---+---+--+-+++-++-+------+-++-+---+--+-+++-+++-+++-++-+--+-++-+--
+++-++++---+---+---+----+++-++++--------++++---+----+++-+++-+++-++++----++++---
--+---++++----+---+---+---++++-+++----+----+++----+---++++-+++-+++-+++----++++-
-+---+-++-+--+---+---+---+-++-+++-+--+-----++-+--+---+-++-+++-+++-+++-+--+-++-+
+---+--+-++-+---+---+---+--+-+++-++-+------+-++-+---+--+-+++-+++-+++-++-+--+-++
++++----++++---+---+---+----+++-++++--------++++---+----+++-+++-+++-++++----+++
--++++----++++----+---+---+---++++-+++----+----+++----+---++++-+++-+++-+++----+
-+-++-+--+-++-+--+---+---+---+-++-+++-+--+-----++-+--+---+-++-+++-+++-+++-+--+-
+--+-++-+--+-++-+---+---+---+--+-+++-++-+------+-++-+---+--+-+++-+++-+++-++-+--
+++-++++----++++---+---+---+----+++-++++--------++++---+----+++-+++-+++-++++---
--+---++++----++++----+---+---+---++++-+++----+----+++----+---++++-+++-+++-+++-
-+---+-++-+--+-++-+--+---+---+---+-++-+++-+--+-----++-+--+---+-++-+++-+++-+++-+
+---+--+-++-+--+-++-+---+---+---+--+-+++-++-+------+-++-+---+--+-+++-+++-+++-++
++++----++++----++++---+---+---+----+++-++++--------++++---+----+++-+++-+++-+++
--++++----++++----++++----+---+---+---++++-+++----+----+++----+---++++-+++-+++-
-+-++-+--+-++-+--+-++-+--+---+---+---+-++-+++-+--+-----++-+--+---+-++-+++-+++-+
+--+-++-+--+-++-+--+-++-+---+---+---+--+-+++-++-+------+-++-+---+--+-+++-+++-++
+++-++++----++++----++++---+---+---+----+++-++++--------++++---+----+++-+++-+++
--++++-+++----++++----++++----+---+---+---++++-+++----+----+++----+---++++-+++-
-+-++-+++-+--+-++-+--+-++-+--+---+---+---+-++-+++-+--+-----++-+--+---+-++-+++-+
+--+-+++-++-+--+-++-+--+-++-+---+---+---+--+-+++-++-+------+-++-+---+--+-+++-++
+++-+++-++++----++++----++++---+---+---+----+++-++++--------++++---+----+++-+++
--++++-+++-+++----++++----++++----+---+---+---++++-+++----+----+++----+---++++-
-+-++-+++-+++-+--+-++-+--+-++-+--+---+---+---+-++-+++-+--+-----++-+--+---+-++-+
+--+-+++-+++-++-+--+-++-+--+-++-+---+---+---+--+-+++-++-+------+-++-+---+--+-++
+++-+++-+++-++++----++++----++++---+---+---+----+++-++++--------++++---+----+++
--++++-+++-+++-+++----++++----++++----+---+---+---++++-+++----+----+++----+---+
-+-++-+++-+++-+++-+--+-++-+--+-++-+--+---+---+---+-++-+++-+--+-----++-+--+---+-
+--+-+++-+++-+++-++-+--+-++-+--+-++-+---+---+---+--+-+++-++-+------+-++-+---+--
+++-+++-+++-+++-++++----++++----++++---+---+---+----+++-++++--------++++---+---
--+---++++-+++-+++-+++----++++----++++----+---+---+---++++-+++----+----+++----+
-+---+-++-+++-+++-+++-+--+-++-+--+-++-+--+---+---+---+-++-+++-+--+-----++-+--+-
+---+--+-+++-+++-+++-++-+--+-++-+--+-++-+---+---+---+--+-+++-++-+------+-++-+--
++++----+++-+++-+++-++++----++++----++++---+---+---+----+++-++++--------++++---
--+---+---++++-+++-+++-+++----++++----++++----+---+---+---++++-+++----+----+++-
-+---+---+-++-+++-+++-+++-+--+-++-+--+-++-+--+---+---+---+-++-+++-+--+-----++-+
+---+---+--+-+++-+++-+++-++-+--+-++-+--+-++-+---+---+---+--+-+++-++-+------+-++
++++---+----+++-+++-+++-++++----++++----++++---+---+---+----+++-++++--------+++
--++++----+---++++-+++-+++-+++----++++----++++----+---+---+---++++-+++----+----
-+-++-+--+---+-++-+++-+++-+++-+--+-++-+--+-++-+--+---+---+---+-++-+++-+--+-----
+--+-++-+---+--+-+++-+++-+++-++-+--+-++-+--+-++-+---+---+---+--+-+++-++-+------
+++-++++---+----+++-+++-+++-++++----++++----++++---+---+---+----+++-++++-------
\end{verbatim}
\end{minipage}
\end{equation*}

\begin{eqnarray*}
\text{Determinant:}&& 2^{78}\times4\times 19^{30}\times 23^9\\
\text{Theoretical bound:}&&2^{78}\times4\times19^{36}\times225\sqrt{40145}\\
\text{Fraction of bound:}&& 0.84924
\end{eqnarray*}
\vfill

Our conference-matrix construction produces this example (see
Sec.~\ref{sec:methods}). We find no prior records in the
literature.

\vfill
\eject

$$\mathbf{n=93:}$$
\begin{equation*}
\begin{minipage}[t]{4.75in}
\tiny
\baselineskip 0pt
\begin{verbatim}
+-----+--+-++---+---+-----+---------+-++--+++--+-------+---+--++----------+-----+-+--+--+----
+-------+-----+--+-----+----+-+++-++-----+----+--+--+----+--+---+----+-----+-++------------++
++++---+-----+-+-----+---+-------+--------------+-++------+--+---++++-+--+--+------+--+------
+---++----+-------++--+-+--+-+-------+--+----+-------++-+--------------++------+-+--+--+-++--
-+----++++-+++--------+-+--+++-++--++-+-+---++----+--+---+---+-+++-++--+--+---+----++------+-
--++++-++--++-+--++------++---+---+--+-+-+-+++--+++--+------+--+----++--+-+----+---++-+-+----
-+-+--+-+-+----+++----++-+-+----++++-+-++-+--++---++-+--+-+++-+-+-----+---+----++-+----------
-+-++++--+--++--+++-+++--++---------++-----------+-+-----++-++-++-----------++++-+---+---++++
--+----+-+-++--+-+--++-+-+++-+-++--++---++------++-+--+++-----+--+-----++----+-+--+---++----+
-++-++--+++--+++-+-+----+-+------+-+++++--+-+----+--+--++--++----++++--+---+-+------------+--
---+--+-+++-+-+++-+-++--+--++-++-++----+---+---+--+--+--------+--+-----+-+-+++------+-+--++-+
-++--+--+-+++-+-+----++-+-+-+--++-----++--+---+-++-+-++-+--+---+---+-++-----+-+-++-+--++-+--+
---++-+++--++-++---+--+----++---------+-++-+--+--+-++--++-+----++++---+----+++-++--+++-++-+--
-+--+++-++-+------++++++--++-++--++------++---++-+------++---++-+++++++---+-----------++++---
-+--+--+-++--+----+++-+-+---+--++++-----++++-+-++++++------++--+---+-----+---+++--+--++++----
-++---+--+-++---+-++-----+--+-+++-+--+--+++-++------+-+++-+-++++--++--+-++-+-------+-----++-+
--+--+-+---+----+--++++-+-+---+-+++++----+--+-+---+++---++++--+--------+-++++-++-+-++---+-+--
-++-+-----+++++--+-+----+---+++-+--+-+----++--++---++----++--++--+--+---+-+-+--++-+-+++--+--+
--+++-+++++-----++-+---+-----+-+---++-+---+-+--++-+-++--+-+--+--+----+---+-+---+-+-+-+++++--+
-++++-+----+----+-+-++-+--+-+---++-+--+----+-+---++-+++--++-+----+++-+-++--+----+++-++-+-----
-+-+---++------+--++--+---+-++++--++++-+---+----+-----++-+-++++--+-+-------++--+++++---+++-+-
-+-+-+-+--++--+-++--+-+--+---+++--+--++-+-+----+-+-----+-+-+----+++---++++-+--+--++++++------
--+--+-++-+-+--+-++-+++----+---+-+---+--+--++-+++---+-+-++++-+-+--++-+----++------+-++-----++
---++-------++-++--++-----+++++++--+-+-+-+----+++-++-+-+++-+---++-++++--++-------+---+----+--
--+--+++-----++-+-+--+-+--+-+++-----+--++-++-++++-+-+--++-++----++-+--+-+----+++----+----+-+-
-+-++-+--+-+-+-+---++-++-+--------++---++--++-++++--+-+-+--+--------++-+++--+-+-+--++----+-++
---+++---+-+-++++------+----++-+++---++-----+-+-+-----+--++++-++--+---++-++--+-+----+-++++-+-
-+--++--------+-++-+++-+++-+-+++-----++----+++----+++-+-+-----++++-+-+---+--++--+-++----+--+-
---+--++-++--++---++-+--++++--+----+--+--+----+-----+-++---+++-++--+-+++-++----++++-+-+-----+
--++++--------+---++++---+-+++--+-+-+-+-+---+-++--++-+-+---+++----+-+-----++-++-+-++---+-+--+
-+--+---+-+--+---++-+----+++++--++--+-++-+-+-++-------+++-+---+-----+-+--+-+--+--+-++++-+--++
--+---++-+-+-++--++---++++-++--+-+---+-+----------+++-++---+--+-+-+-++--+-+++-+--+---+-+++---
--++-+--++---++-+--+---+--------++---+--++-+-+++++---+++-----++-++-+---+--+-+-+--+++-+--+-+-+
---+-+-+-++++----+-+-++++----+--+---+--+-+-+-+-+-----++--++++---+-+++-+-++-+++--+-------+---+
-----++-+-+-+++----+--++--+-+++--++-+---+---+---++++--+---+--+++----++++++-+-+--+++--+-------
-++----+---+-+--++------+-++--++-++---++++-++-++-+---++---+--+--+-+-+----+---+-++++----+-+++-
-+---+++------++-+-+-+-+++--+---+-+-+--+-----+-+++--++-+++---+++--+------++-----++--++++--++-
----++++--+-++-----+++---------+++++-+++-+--++--+------+--+-----++---++-+-+++-+-+-+---++-+++-
-++--++---+---++--+-+--+++-++--+---++----++++-+++--------++-+--++---+-+++---+----+++---++-+--
----+--+--+-+++++---+-+++++---+-++--++--++++------+-++-+-+---+----+--+++--------+--+---+++-++
--+-+-+++-+-++-++-+--+--+----+----++--+-+-+--+++-+-+--++-+--+-+---+--+--+++-++---+-+----+--+-
--+-+----+-+--++-++--++----++-+---+++++--++--+-++--+-++---++-+---+---+++-+----+-+----+--+-++-
-+-+-+--+--+---+--+----++----++--+-++-++++--++-+--+++-++-+-----+--+--++-----+-++--+--++--++-+
--++----+-++-+-+---+++++-+---+++------++-++++------+++------+++----+--+-+-+--++--+--++-+--++-
--++-----++-+---+------+++++-++--++---+-+----++-+---+---++-+++-+-+--+---+--++-+-+----++++-++-
----+-++-------++++---+++-+--+-+-+--+++--+----++++-+++------+-++-+-++-+-++++--+-+-+++----++-+
--++--+-++----+---+-+-+-+-+--+--+----+-++++-++--+--++----+++--+-+++-++---+---+--+++++-++--+++
-+--------++++---++-+--+--+-++-+-++---------++-++--+++-++--+++--+++--+++--+-++-+++-++---+-+-+
--+++--++--+-+------++-+++--+--+--+-++-++-+--++---------++-----+-+-+++++-+++-+-++++--+--+++-+
----+--+++----++++--++------++---++----++++-+------++++++++-++-++--++--++---+-+-+++-+-+-++---
---+--+----+-+++-+---+--+-++-+--++++-+---++-++-+-++-----+--+-+-+---+--+-+--+++---+++++-+++-++
-+++---+--+-----+----++----++++--+-+++---+++++--++-----++---+-+++-+-++-+-++-+++-+---++-+-+--+
-----------+++++--++-+---+++-+-+-+--++-++-+-++++-+--++---++-+-+-+--+-+-+-+-++-+++--+-+++-----
--+--++-----++-+--+---++-+----+-++++-++---+++--+---+-+++++-----+++-++----+-+-+++++--+-+++---+
-+-+-----+--+-----++-+--+-----++-++-++++--+---+-+++++++++++-----+---+-+++-+--+----++++-++--++
-+----+-+---+++--+---++-++---+----+---+--+-++-+-+-+--+-++++++-+--+++-+-+++----++--+--+-++++-+
-----+---+++-+-++----++-----+-----+-++++-+-+-+++--+++---+---+++++++-++-++--+---+-+++--++---++
---+--------++--++++--++++-++-+-+---+-+++++-+--+++-+----+-+-+----+++-+-+--+++----++-+-+-++-+-
--+-++--+--++--++--+--+-++++---+--++---------+-++-+-+--+--+-+-+-+-+++-++-----++-+++-+++-++-++
----+---++-------++--++++-+---+++--+-+-++--++-++--+----++-+-++-+-++---+++++-+++-++-+-+++-----
--++-+---+--+--+-+--+-+-+--+++--+++-+-----+++-+--++-+-++----+---++-+-++-+++----+++-+-++--+++-
-++++----++--++-------++--++---+--+-+----++----+--+--++++++---++--++-+--+-+++++++-+++-+---++-
---+-+-+--+---+---+-+-++--++---++-++--++--+-+++--+-+++-+--+--+++-++++--+++--+--+-----+--++-++
-----+-+-+++-+-+--+-+---++---++-+--+---+--+---+-+++--++-+-+-+++++++--+-----+-++-+-++++--+++--
-+---+-+++---+--+---++-++--++-+----+-+--+---+--+-+-+-+-+--+++-+----++++-+--+---++--++++++-+++
-+-+++-++---+--++-+-----+++-+---+--------++-+--+----+++-+--+-++-++--++++++++++++--+--+-+---+-
---+---++---+----+-++--++++-+----++++++---++-+++---+--+---++-++++-++-+-++----++----++-++-++--
--+-+-----++----++++--+--+--+----+-++--++---+-+--++--+-+-+-+-+++++-+-++-++--++----+-+-+++-+++
-+---+---+--+--+-----+++-++-++-+---+-++-++++-++--++-++-+--++-+------++---+++++---++-+-+-+++--
-++------++-+-++++++-+-+--+--+--+-++---++--+-------+-----------++-++++++-+++--++--++++-+++++-
---+-++-+-++-+-++++--------+--+++---+------+-+---+-++--+++++-+-----+++-++++--++-+-++--+++++--
-----++---++---+++-++--++-+++-----+--++++++-----+-+---++-++-++---+-+++---++-++++---+++-+----+
--+---++-++-+--+-++++--+----+-+------++--+--+++--++--+--++-+--+---+++-++----++++++++-++--+-+-
---++++---+-+----+---+---+--++++-+-+---++++----+-++-+-+---++--++-++++--+--+---++++-+---++-+++
--+---+-------+-+--++-+-++---+-+-+-+---+-+-+-+-+-+-+-+++++++++-+-+--++++--++-+-+---+-++--+-+-
-------++--++++-+-+++--+---+-+----++++----++-----++-+++-+-+++--+-++-+-+--++-+-+--++---++--+++
-+------+++-+-+-+---+----+-+-------+++-+-+--+++++++++-+--+---+---+++--+-++++-+++++--+--+++---
---+---------++++++++++-+-+-+--+---+--+-+-++++--+++---+-++---++-+-+-+-+-+-++--+++-----+-+-+-+
-+--+-++----+-+++--+-+++---+--++++--+--+--+-+++-------+----+++---++-++--+----+++-++++++-+-+-+
-+--++-+---++-+---+----+---+-----+-+--++++++----+-++----++-++++-+--+---+++++-+--++-+-++--++++
----++++-+---+--++----+--++++++---++---+--+++++-+--++++--+--------++--++--++++----++-+++--+-+
--+--+++++---------+-----+++++-+++---+++--++---+-+++----++--++----+--+++++--+-+++-+-+---+-+++
-----+--++---+--++++-+---+-+---++++++-+-++++----+---++---+-+--++-++-+++-+---++-+++-+++---+---
-++-+-+--------+----+---+--+-++-+-+-+++++--+----++--++----++-++++-++-++++-++++++-+----+-+---+
-+-----+++-+--+++++++-+--+---++-++----+---+--++++-+----+--+----+---+++-++-++++-+++-----+--+++
--+-+---+---++---+--++++---+--+-+-----++--+--+-++++-+-++-+-++++++---+-++-+-++--+--+++--+-++--
---++++--+-++-+--+--+---+-----++-+-++---+-++-++-+--+-+-++---+++---+++++--+-++-++-++-+--+---+-
----+-+-+-++----+-+--+-+++---+-++-+-++---++++-+-++++--++---+-+--+-+++--+---++--++---+++---++-
----+--+-++---++---------++---+++-+-+-+-+--+++-+-+++--+-+++++-+-++-++++---+-+----+--++-+-++-+
-++--++-+----+-+-+-++-----+---++----+---+++--++---++-++----++--++++----++++++---+---++++++-++
------+-+-++--+-----+++--++-+-+-++-++++-+------++---+++-+---+---+-+++---+++-++-+--++-++-++++-
-+-+--++--++-----+-+-+----++----+----+--+-+++-+++-++++++-+--+-++-+--++++---+-++--+----+-++++-
----+---++++---+---+-+-++++++-+--+--------++++--++-+-+-+-+++---++++-----++++--++-+++-----+-++
\end{verbatim}
\end{minipage}
\end{equation*}
\begin{eqnarray*}
\text{Determinant:}&& 2^{92}\times23^{45}\times 231\\
\text{Theoretical bound:}&&2^{92}\times 23^{46}\times\sqrt{185}\\
\text{Fraction of bound:}&& 0.738411
\end{eqnarray*}
\vfill
Farmakis \& Kounias~\cite{fk} published an earlier record of
$\,2^{92}\times23^{45}\times230\,$ (73.52\% of bound).
\eject
$$\mathbf{n=95:}$$
\begin{equation*}
\begin{minipage}[t]{4.8in}
\tiny
\baselineskip 0pt
\begin{verbatim}
+++--++--++--++--++--++--++--++--++--++--++--++--++--++--++--++--++--++--++--++--++--++--++--++
+++-+-+-+-+-+-+-+-+-+-+-+-+-+-+-+-+-+-+-+-+-+-+-+-+-+-+-+-+-+-+-+-+-+-+-+-+-+-+-+-+-+-+-+-+-+-+
++++--++--++--++--++--++--++--++--++--++--++--++--++--++--++--++--++--++--++--++--++--++--++--+
--+----+++-+++-+++-+++----++++----++++-+++----+---++++-+++----+---++++----++++----+---+---+---+
-+-----++-+++-+++-+++-+--+-++-+--+-++-+++-+--+---+-++-+++-+--+---+-++-+--+-++-+--+---+---+---+-
+------+-+++-+++-+++-++-+--+-++-+--+-+++-++-+---+--+-+++-++-+---+--+-++-+--+-++-+---+---+---+--
+++-----+++-+++-+++-++++----++++----+++-++++---+----+++-++++---+----++++----++++---+---+---+---
--+---+----+++-+++-+++-+++----++++----++++-+++----+---++++-+++----+---++++----++++----+---+---+
-+---+-----++-+++-+++-+++-+--+-++-+--+-++-+++-+--+---+-++-+++-+--+---+-++-+--+-++-+--+---+---+-
+---+------+-+++-+++-+++-++-+--+-++-+--+-+++-++-+---+--+-+++-++-+---+--+-++-+--+-++-+---+---+--
++++--------+++-+++-+++-++++----++++----+++-++++---+----+++-++++---+----++++----++++---+---+---
--+---+---+----+++-+++-+++-+++----++++----++++-+++----+---++++-+++----+---++++----++++----+---+
-+---+---+-----++-+++-+++-+++-+--+-++-+--+-++-+++-+--+---+-++-+++-+--+---+-++-+--+-++-+--+---+-
+---+---+------+-+++-+++-+++-++-+--+-++-+--+-+++-++-+---+--+-+++-++-+---+--+-++-+--+-++-+---+--
++++---+--------+++-+++-+++-++++----++++----+++-++++---+----+++-++++---+----++++----++++---+---
--+---+---+---+----+++-+++-+++-+++----++++----++++-+++----+---++++-+++----+---++++----++++----+
-+---+---+---+-----++-+++-+++-+++-+--+-++-+--+-++-+++-+--+---+-++-+++-+--+---+-++-+--+-++-+--+-
+---+---+---+------+-+++-+++-+++-++-+--+-++-+--+-+++-++-+---+--+-+++-++-+---+--+-++-+--+-++-+--
++++---+---+--------+++-+++-+++-++++----++++----+++-++++---+----+++-++++---+----++++----++++---
--+---+---+---+---+----+++-+++-+++-+++----++++----++++-+++----+---++++-+++----+---++++----++++-
-+---+---+---+---+-----++-+++-+++-+++-+--+-++-+--+-++-+++-+--+---+-++-+++-+--+---+-++-+--+-++-+
+---+---+---+---+------+-+++-+++-+++-++-+--+-++-+--+-+++-++-+---+--+-+++-++-+---+--+-++-+--+-++
++++---+---+---+--------+++-+++-+++-++++----++++----+++-++++---+----+++-++++---+----++++----+++
--++++----+---+---+---+----+++-+++-+++-+++----++++----++++-+++----+---++++-+++----+---++++----+
-+-++-+--+---+---+---+-----++-+++-+++-+++-+--+-++-+--+-++-+++-+--+---+-++-+++-+--+---+-++-+--+-
+--+-++-+---+---+---+------+-+++-+++-+++-++-+--+-++-+--+-+++-++-+---+--+-+++-++-+---+--+-++-+--
+++-++++---+---+---+--------+++-+++-+++-++++----++++----+++-++++---+----+++-++++---+----++++---
--+---++++----+---+---+---+----+++-+++-+++-+++----++++----++++-+++----+---++++-+++----+---++++-
-+---+-++-+--+---+---+---+-----++-+++-+++-+++-+--+-++-+--+-++-+++-+--+---+-++-+++-+--+---+-++-+
+---+--+-++-+---+---+---+------+-+++-+++-+++-++-+--+-++-+--+-+++-++-+---+--+-+++-++-+---+--+-++
++++----++++---+---+---+--------+++-+++-+++-++++----++++----+++-++++---+----+++-++++---+----+++
--++++----++++----+---+---+---+----+++-+++-+++-+++----++++----++++-+++----+---++++-+++----+---+
-+-++-+--+-++-+--+---+---+---+-----++-+++-+++-+++-+--+-++-+--+-++-+++-+--+---+-++-+++-+--+---+-
+--+-++-+--+-++-+---+---+---+------+-+++-+++-+++-++-+--+-++-+--+-+++-++-+---+--+-+++-++-+---+--
+++-++++----++++---+---+---+--------+++-+++-+++-++++----++++----+++-++++---+----+++-++++---+---
--+---++++----++++----+---+---+---+----+++-+++-+++-+++----++++----++++-+++----+---++++-+++----+
-+---+-++-+--+-++-+--+---+---+---+-----++-+++-+++-+++-+--+-++-+--+-++-+++-+--+---+-++-+++-+--+-
+---+--+-++-+--+-++-+---+---+---+------+-+++-+++-+++-++-+--+-++-+--+-+++-++-+---+--+-+++-++-+--
++++----++++----++++---+---+---+--------+++-+++-+++-++++----++++----+++-++++---+----+++-++++---
--+---+---++++----++++----+---+---+---+----+++-+++-+++-+++----++++----++++-+++----+---++++-+++-
-+---+---+-++-+--+-++-+--+---+---+---+-----++-+++-+++-+++-+--+-++-+--+-++-+++-+--+---+-++-+++-+
+---+---+--+-++-+--+-++-+---+---+---+------+-+++-+++-+++-++-+--+-++-+--+-+++-++-+---+--+-+++-++
++++---+----++++----++++---+---+---+--------+++-+++-+++-++++----++++----+++-++++---+----+++-+++
--++++----+---++++----++++----+---+---+---+----+++-+++-+++-+++----++++----++++-+++----+---++++-
-+-++-+--+---+-++-+--+-++-+--+---+---+---+-----++-+++-+++-+++-+--+-++-+--+-++-+++-+--+---+-++-+
+--+-++-+---+--+-++-+--+-++-+---+---+---+------+-+++-+++-+++-++-+--+-++-+--+-+++-++-+---+--+-++
+++-++++---+----++++----++++---+---+---+--------+++-+++-+++-++++----++++----+++-++++---+----+++
--++++-+++----+---++++----++++----+---+---+---+----+++-+++-+++-+++----++++----++++-+++----+---+
-+-++-+++-+--+---+-++-+--+-++-+--+---+---+---+-----++-+++-+++-+++-+--+-++-+--+-++-+++-+--+---+-
+--+-+++-++-+---+--+-++-+--+-++-+---+---+---+------+-+++-+++-+++-++-+--+-++-+--+-+++-++-+---+--
+++-+++-++++---+----++++----++++---+---+---+--------+++-+++-+++-++++----++++----+++-++++---+---
--+---++++-+++----+---++++----++++----+---+---+---+----+++-+++-+++-+++----++++----++++-+++----+
-+---+-++-+++-+--+---+-++-+--+-++-+--+---+---+---+-----++-+++-+++-+++-+--+-++-+--+-++-+++-+--+-
+---+--+-+++-++-+---+--+-++-+--+-++-+---+---+---+------+-+++-+++-+++-++-+--+-++-+--+-+++-++-+--
++++----+++-++++---+----++++----++++---+---+---+--------+++-+++-+++-++++----++++----+++-++++---
--+---+---++++-+++----+---++++----++++----+---+---+---+----+++-+++-+++-+++----++++----++++-+++-
-+---+---+-++-+++-+--+---+-++-+--+-++-+--+---+---+---+-----++-+++-+++-+++-+--+-++-+--+-++-+++-+
+---+---+--+-+++-++-+---+--+-++-+--+-++-+---+---+---+------+-+++-+++-+++-++-+--+-++-+--+-+++-++
++++---+----+++-++++---+----++++----++++---+---+---+--------+++-+++-+++-++++----++++----+++-+++
--++++----+---++++-+++----+---++++----++++----+---+---+---+----+++-+++-+++-+++----++++----++++-
-+-++-+--+---+-++-+++-+--+---+-++-+--+-++-+--+---+---+---+-----++-+++-+++-+++-+--+-++-+--+-++-+
+--+-++-+---+--+-+++-++-+---+--+-++-+--+-++-+---+---+---+------+-+++-+++-+++-++-+--+-++-+--+-++
+++-++++---+----+++-++++---+----++++----++++---+---+---+--------+++-+++-+++-++++----++++----+++
--++++-+++----+---++++-+++----+---++++----++++----+---+---+---+----+++-+++-+++-+++----++++----+
-+-++-+++-+--+---+-++-+++-+--+---+-++-+--+-++-+--+---+---+---+-----++-+++-+++-+++-+--+-++-+--+-
+--+-+++-++-+---+--+-+++-++-+---+--+-++-+--+-++-+---+---+---+------+-+++-+++-+++-++-+--+-++-+--
+++-+++-++++---+----+++-++++---+----++++----++++---+---+---+--------+++-+++-+++-++++----++++---
--+---++++-+++----+---++++-+++----+---++++----++++----+---+---+---+----+++-+++-+++-+++----++++-
-+---+-++-+++-+--+---+-++-+++-+--+---+-++-+--+-++-+--+---+---+---+-----++-+++-+++-+++-+--+-++-+
+---+--+-+++-++-+---+--+-+++-++-+---+--+-++-+--+-++-+---+---+---+------+-+++-+++-+++-++-+--+-++
++++----+++-++++---+----+++-++++---+----++++----++++---+---+---+--------+++-+++-+++-++++----+++
--++++----++++-+++----+---++++-+++----+---++++----++++----+---+---+---+----+++-+++-+++-+++----+
-+-++-+--+-++-+++-+--+---+-++-+++-+--+---+-++-+--+-++-+--+---+---+---+-----++-+++-+++-+++-+--+-
+--+-++-+--+-+++-++-+---+--+-+++-++-+---+--+-++-+--+-++-+---+---+---+------+-+++-+++-+++-++-+--
+++-++++----+++-++++---+----+++-++++---+----++++----++++---+---+---+--------+++-+++-+++-++++---
--+---++++----++++-+++----+---++++-+++----+---++++----++++----+---+---+---+----+++-+++-+++-+++-
-+---+-++-+--+-++-+++-+--+---+-++-+++-+--+---+-++-+--+-++-+--+---+---+---+-----++-+++-+++-+++-+
+---+--+-++-+--+-+++-++-+---+--+-+++-++-+---+--+-++-+--+-++-+---+---+---+------+-+++-+++-+++-++
++++----++++----+++-++++---+----+++-++++---+----++++----++++---+---+---+--------+++-+++-+++-+++
--++++----++++----++++-+++----+---++++-+++----+---++++----++++----+---+---+---+----+++-+++-+++-
-+-++-+--+-++-+--+-++-+++-+--+---+-++-+++-+--+---+-++-+--+-++-+--+---+---+---+-----++-+++-+++-+
+--+-++-+--+-++-+--+-+++-++-+---+--+-+++-++-+---+--+-++-+--+-++-+---+---+---+------+-+++-+++-++
+++-++++----++++----+++-++++---+----+++-++++---+----++++----++++---+---+---+--------+++-+++-+++
--++++-+++----++++----++++-+++----+---++++-+++----+---++++----++++----+---+---+---+----+++-+++-
-+-++-+++-+--+-++-+--+-++-+++-+--+---+-++-+++-+--+---+-++-+--+-++-+--+---+---+---+-----++-+++-+
+--+-+++-++-+--+-++-+--+-+++-++-+---+--+-+++-++-+---+--+-++-+--+-++-+---+---+---+------+-+++-++
+++-+++-++++----++++----+++-++++---+----+++-++++---+----++++----++++---+---+---+--------+++-+++
--++++-+++-+++----++++----++++-+++----+---++++-+++----+---++++----++++----+---+---+---+----+++-
-+-++-+++-+++-+--+-++-+--+-++-+++-+--+---+-++-+++-+--+---+-++-+--+-++-+--+---+---+---+-----++-+
+--+-+++-+++-++-+--+-++-+--+-+++-++-+---+--+-+++-++-+---+--+-++-+--+-++-+---+---+---+------+-++
+++-+++-+++-++++----++++----+++-++++---+----+++-++++---+----++++----++++---+---+---+--------+++
--++++-+++-+++-+++----++++----++++-+++----+---++++-+++----+---++++----++++----+---+---+---+----
-+-++-+++-+++-+++-+--+-++-+--+-++-+++-+--+---+-++-+++-+--+---+-++-+--+-++-+--+---+---+---+-----
+--+-+++-+++-+++-++-+--+-++-+--+-+++-++-+---+--+-+++-++-+---+--+-++-+--+-++-+---+---+---+------
+++-+++-+++-+++-++++----++++----+++-++++---+----+++-++++---+----++++----++++---+---+---+-------
\end{verbatim}
\end{minipage}
\end{equation*}
\begin{eqnarray*}
\text{Determinant:}&& 2^{96}\times 23^{36}\times27^{11}\\
\text{Theoretical bound:}&&2^{96}\times 23^{44}\times333\sqrt{69153}\\
\text{Fraction of bound:}&& 0.810642
\end{eqnarray*}
Our conference-matrix construction produces this example (see
Sec.~\ref{sec:methods}). We find no prior records in the
literature.
\vfill
\goodbreak